\documentclass[11pt]{article}

\usepackage[utf8]{inputenc}
\usepackage{fullpage}
\usepackage{amsmath,amsthm,amsfonts,amssymb,amstext,mathrsfs,enumerate}
\usepackage{mathtools}
\usepackage{graphicx,ragged2e,framed,xcolor}

\usepackage{tikz}
\usetikzlibrary{arrows,decorations.markings,shadows}
\usepackage{array,booktabs,longtable,pdflscape}
\usepackage{float}
\floatstyle{ruled}
\newfloat{algorithm}{tbp}{loa}
\floatname{algorithm}{Algorithm}
\usepackage[short]{optidef}
\usepackage{xurl}
\usepackage{authblk}

\theoremstyle{plain}
\newtheorem{theorem}{Theorem}[section]

\newtheorem{proposition}[theorem]{Proposition}

\theoremstyle{definition}

\numberwithin{equation}{section}
\allowdisplaybreaks

\usepackage[pagebackref,hypertexnames=false]{hyperref}
\hypersetup{
  colorlinks=true,
  urlcolor=purple,
  linkcolor=purple,
  citecolor=purple,
  pdftitle={Branch-and-Cut for Maximum k-Independent Sets with Spectral Bounds},
  pdfauthor={Aida Abiad, Julian Golak, Arne Schulz}
}

\DeclareMathOperator{\dist}{dist}
\newcommand{\alphaopt}[1]{\ensuremath{#1^{\ast}}}
\newcommand{\alphanon}[1]{\ensuremath{#1}}
\newcommand{\affl}[3]{\noindent #1, Email: {\ttfamily #2}\\ \textsc{#3}\par\vspace{1.5pt}}

\title{\textbf{Branch-and-Cut for Maximum $k$-Independent Sets with Spectral Bounds}}
\author[1,2]{Aida Abiad}
\author[3]{Julian Golak}
\author[3]{Arne Schulz}
\affil[1]{Department of Mathematics and Computer Science, Eindhoven University of Technology, Eindhoven, The Netherlands}
\affil[2]{Department of Mathematics and Data Science, Vrije Universiteit Brussel, Brussel, Belgium}
\affil[3]{Institute of Operations Management, University of Hamburg Business School, Hamburg, Germany}
\date{}

\begin{document}
\maketitle

\begin{abstract}
We develop a branch-and-cut method for finding maximum distance-independent sets in graphs.
Our approach solves a sequence of fixed-size clique problems, combining graph reductions with combinatorial and spectral bounds.
We give a compact formulation, derive strengthened inequalities, establish conditions for using spectral information, and evaluate them in a comprehensive computational study.
\end{abstract}



\section{Introduction}\label{sec:introduction}
Let $G$ be a simple undirected graph and let $k$ be a positive integer. A set of vertices is \emph{$k$-independent} if the distance between any two distinct vertices in the set exceeds $k$. Its maximum possible cardinality is the \emph{$k$-independence number}, denoted by $\alpha_k(G)$. This generalizes the usual independence number and describes separation requirements in coding theory and location problems. In particular, distance-independent sets in Hamming and sum-rank-metric graphs correspond to codes with a prescribed minimum distance~\cite{abiad2023eigenvalue}.

Finding a maximum $k$-independent set is NP-hard for every fixed positive $k$~\cite{Kong1993OnSets,kz2000,AK2023,Eto2014}. The problem can be expressed as an independent set problem in the graph power $G^k$, which joins vertices whose distance in $G$ is at most $k$, or as a clique problem in its complement~\cite{abiad2019kindependence,ACFNS2022}. This permits the use of established maximum clique algorithms; see the surveys~\cite{bomze1999maximum,wu2015review}. Integer programming provides another approach, with the advantage that additional inequalities and side constraints can be incorporated directly. Martins~\cite{martins2010extended} develops the size-indexed formulation used as our computational reference. Sequential fixed-size searches also appear in methods for related graph optimization problems~\cite{veremyev2025ultrasmall}.

Spectral graph theory provides a complementary source of information. Bounds for the $k$-independence number derived from the spectrum of the original graph have been developed in~\cite{abiad2016spectral,abiad2019kindependence,ACFNS2022,abiad2023eigenvalue,abiad2024unified,abiad2026inertia}. Such bounds may rule out larger solutions before an expensive search is required. They can also be applied locally to restrict how many vertices a solution may select from a given subset. This second use connects spectral bounds with the classical rank inequalities and lifting methods used in stable set optimization~\cite{chvatal1975rank,nemhauser1974properties,rossi2001branchcut,rebennack2011branchcut}.

In this paper, we develop a branch-and-cut framework that combines these ideas for distance independence. We solve a sequence of fixed-size clique problems on reduced auxiliary graphs. Compared with Martins' formulation, each subproblem avoids maintaining separate copies of the vertex variables for all candidate clique sizes. We establish the correctness of the formulation and derive stronger combinatorial and spectral inequalities tailored to the chosen solution size. We also show how information about a subset of vertices can restrict the choices made elsewhere in the graph. The contribution lies in this formulation and its inequalities, together with their integration into the sequential method.

The spectral inequalities use established independence bounds on the power graph and on complements of subgraphs of the current working graph. This differs from optimizing polynomial bounds directly from the spectrum of the original graph. We show how to use the local bounds after vertex and edge reductions and how to strengthen the resulting inequalities. These results strengthen the optimization model while using established eigenvalue bounds.

We evaluate a reference implementation on several Sagemath named graphs and  coding-theoretic graphs such as sum-rank-metric graphs for four distance thresholds. The configuration combining both cut families reports optimality in 470 of 544 cases, compared with 416 for the tested discretized formulation. We examine the contribution of the cut families and the time spent on successive target sizes, and include additional comparisons on DIMACS graphs.%

\section{Reduction to the maximum clique problem}\label{sec:MCP}

\paragraph{Notation}
All graphs are finite, simple, and undirected. We assume $|V|\geq1$ unless stated otherwise; for the graph with no vertices, $\alpha_k(G)=0$.
The distance $\dist_G(i,j)$ is the number of edges on a shortest path from $i$ to $j$, and is infinite when the vertices lie in different connected components.
When such a path exists, $SP_G(i,j)$ denotes one shortest path, viewed as an ordered vertex sequence.
For a graph $H=(U,F)$, write $H[W]$ for the subgraph induced by $W\subseteq U$, $\overline H$ for its complement, $N_H(i)$ for the open neighborhood of $i$, and $d_H(i)=|N_H(i)|$ for its degree.
We write $y(W)=\sum_{i\in W}y_i$ and use $\omega(H)$ and $\alpha_1(H)$ for the clique and independence numbers.
The notation $\bar G=\overline{G^k}$ is reserved for the original auxiliary graph; $H$ denotes a current working graph and $H_q=(V_q,E_q)$ the graph after reductions for target size $q$.

For a graph $H$, let $\omega(H)$ denote its clique number, with $\omega(\varnothing)=0$.
Given an instance $G=(V,E)$ of \textsc{$k$-MIS}, define the auxiliary graph
\[
\bar G=(V,\bar E),\qquad
\bar E=\bigl\{\{i,j\}\in\tbinom V2:\dist_G(i,j)>k\bigr\}.
\]
Thus $\bar G=\overline{G^k}$, and the definitions give
\begin{equation}\label{eq:power-clique-equivalence}
\alpha_k(G)=\alpha_1(G^k)=\omega(\bar G).
\end{equation}
Indeed, a set of vertices has pairwise distance greater than $k$ in $G$ if and only if all its pairs are edges of $\bar G$~\cite{ACFNS2022,abiad2019kindependence}.
The auxiliary graph can be constructed in polynomial time by computing shortest-path distances~\cite{korte2008combinatorial}; vertices in different connected components of $G$ are adjacent in $\bar G$.
The empty input graph has $\alpha_k(G)=0$ and is handled directly.
For the remainder of this section, assume $V\neq\varnothing$.

\subsection{A compact formulation for a fixed clique size}

The discretized formulation of Martins~\cite{martins2010extended} represents the possible clique sizes explicitly.
Let $Q\subseteq\{1,\ldots,|V|\}$ contain $\omega(\bar G)$.
The binary variable $w_q$ selects size $q$, and $x_{iq}$ indicates whether vertex $i$ belongs to the clique selected at that size.
The formulation is
\begin{maxi!}
{} {\sum_{i\in V}\sum_{q\in Q}x_{iq}}
{}{\label{obj}}
\addConstraint{\sum_{q\in Q}w_q}{=1}\label{Con1}
\addConstraint{\sum_{j\in N_{\bar G}(i)}x_{jq}}{\geq(q-1)x_{iq}\qquad}{\forall i\in V,\ q\in Q}\label{Con2}
\addConstraint{\sum_{i\in V}x_{iq}}{=qw_q}{\forall q\in Q}\label{Con3}
\addConstraint{w_q}{\in\{0,1\}}{\forall q\in Q}\label{Con4}
\addConstraint{x_{iq}}{\in\{0,1\}}{\forall i\in V,\ q\in Q}\label{Con5}
\end{maxi!}
Exactly one size is active.
At that size, every selected vertex has the other $q-1$ selected vertices as neighbors, so the selected set is a clique.
The assumption on $Q$ therefore makes the optimal objective value equal to $\omega(\bar G)$.
The number of variables grows with $|Q|$, which motivates solving one size at a time.

Let $H=(U,F)$ be the working auxiliary graph: initially $H=\bar G$, and after preprocessing it is the reduced graph defined in Section~\ref{sec:rgs}.
We write $N_H(i)$ and $d_H(i)$ for its neighborhoods and degrees.
For a fixed integer $q\geq1$, consider
\begin{maxi!}
{}{\sum_{i\in U}y_i}
{}{\label{objN}}
\addConstraint{\sum_{j\in N_H(i)}y_j}{\geq(q-1)y_i\qquad}{\forall i\in U}\label{ConN1}
\addConstraint{\sum_{i\in U}y_i}{=qz}\label{ConN2}
\addConstraint{y_i}{\in\{0,1\}}{\forall i\in U}\label{ConN4}
\addConstraint{z}{\in[0,1].}{}\label{ConN5}
\end{maxi!}
Here $y$ selects the vertices, while $z$ permits either a $q$-clique or the empty solution.
The variable $z$ is continuous.

\begin{proposition}[Exactness of the fixed-size formulation]\label{prop:fixed-size-exactness}
For every finite simple graph $H$ and every integer $q\geq1$, the feasible solutions of~\eqref{objN}--\eqref{ConN5} are precisely $(y,z)=(0,0)$ and $(y,z)=(\mathbf{1}_C,1)$, where $C$ is a clique of size $q$ in $H$.
Consequently, its optimal objective value is $q$ if such a clique exists and $0$ otherwise.
\end{proposition}
\begin{proof}
Let $(y,z)$ be feasible and set $C=\{i\in U:y_i=1\}$.
If $C=\varnothing$, then~\eqref{ConN2} and $q\geq1$ give $z=0$.
Otherwise, choose $i\in C$.
Constraint~\eqref{ConN1} gives at least $q-1$ selected neighbors of $i$, and hence $|C|\geq q$.
On the other hand,~\eqref{ConN2} and $z\leq1$ give $|C|\leq q$.
Thus $|C|=q$ and $z=1$.
Applying~\eqref{ConN1} to each $i\in C$ shows that $i$ is adjacent to every vertex of $C\setminus\{i\}$, so $C$ is a clique.
Conversely, the empty solution and the incidence vector of any $q$-clique satisfy every constraint.
The argument also covers $q=1$ and $U=\varnothing$.
\end{proof}

The formulation has $|U|$ binary variables, one continuous variable, and $|U|+1$ linear constraints apart from variable bounds, with $O(|U|+|F|)$ nonzero coefficients.
Eliminating $z$ gives the equivalent formulation obtained by replacing~\eqref{ConN2} and~\eqref{ConN5} with $\sum_{i\in U}y_i\leq q$.
We retain $z$ because it provides a useful strengthening of the cuts in Section~\ref{sec:vi}.
The LP relaxation also implies
\begin{equation}\label{eq:y-bounded-by-z}
0\leq y_i\leq z\qquad(i\in U):
\end{equation}
indeed,~\eqref{ConN1} and~\eqref{ConN2} give
$(q-1)y_i\leq\sum_{j\in N_H(i)}y_j\leq qz-y_i$.

\paragraph{Relation to the experimental implementation}
The implementation used in Section~\ref{sec:IPkMIS} declares $z$ binary and retains the additional rows
\begin{equation}\label{ConN3}
\sum_{j\in N_H(i)}y_j\leq q-y_i\qquad(i\in U).
\end{equation}
Proposition~\ref{prop:fixed-size-exactness} shows that the integrality restriction on $z$ does not change the integer feasible solutions.
Moreover,~\eqref{ConN3} is redundant even in the LP relaxation, since
\[
\sum_{j\in N_H(i)}y_j
\leq\sum_{j\in U\setminus\{i\}}y_j
=qz-y_i\leq q-y_i.
\]
The compact formulation and the implementation therefore have the same integer feasible solutions and the same LP relaxation.

\subsection{Reduction of graph size}\label{sec:rgs}

For $q\geq2$, two deletion rules reduce the auxiliary graph before the fixed-size problem is solved:
a vertex $i$ may be deleted when $d_H(i)<q-1$, and an edge $\{i,j\}$ may be deleted when $|N_H(i)\cap N_H(j)|<q-2$.
Degrees and common neighborhoods are evaluated in the current graph, so a deletion may enable further deletions.
Algorithm~\ref{Alg:GR} applies the rules until neither is possible.
Related degree and common-neighbor restrictions appear in Martins~\cite{martins2010extended}; core and truss reductions are standard preprocessing tools for maximum clique~\cite{cohen2008trusses,verma2015solving}.
The case $q=1$ needs no reduction: a nonempty graph contains a one-vertex clique.

\begin{algorithm}[t]
\caption{Graph reduction for a target clique size}\label{Alg:GR}
\begin{tabular}{@{}r@{\quad}p{0.93\linewidth}@{}}
1 & \textbf{Input:} auxiliary graph $\bar G=(V,\bar E)$ and integer $q\geq2$ \\
2 & $H=(U,F)\gets\bar G$ \\
3 & \textbf{while} a vertex or edge of $H$ satisfies a deletion rule \textbf{do} \\
4 & \hspace*{1em}\textbf{if} there is a vertex $i\in U$ with $d_H(i)<q-1$ \textbf{then} \\
5 & \hspace*{2em}delete $i$ and its incident edges from $H$ \\
6 & \hspace*{1em}\textbf{else} \\
7 & \hspace*{2em}choose $\{i,j\}\in F$ with $|N_H(i)\cap N_H(j)|<q-2$ \\
8 & \hspace*{2em}delete $\{i,j\}$ from $H$ \\
9 & \hspace*{1em}\textbf{end if} \\
10 & \textbf{end while} \\
11 & \textbf{Output:} $H_q=(V_q,E_q)\gets H$
\end{tabular}
\end{algorithm}

\begin{proposition}[Validity and uniqueness of the reduction]\label{prop:peeling}
For $q\geq2$, Algorithm~\ref{Alg:GR} terminates and preserves every clique of $\bar G$ of size at least $q$.
Its output is independent of the order of valid deletions and is the unique largest subgraph of $\bar G$ in which every vertex has degree at least $q-1$ and the endpoints of every edge have at least $q-2$ common neighbors.
\end{proposition}
\begin{proof}
Each iteration deletes a vertex or an edge, so the algorithm terminates after at most $|V|+|\bar E|$ iterations.
Let $C$ be a clique of size at least $q$.
As long as all its vertices and edges remain, each vertex of $C$ has at least $q-1$ neighbors in $C$, and each edge of $C$ has at least $q-2$ common neighbors in $C$.
Neither rule can therefore remove the first vertex or edge of $C$.

For uniqueness, let $J$ be any subgraph of $\bar G$ satisfying both degree conditions in the statement.
If $J$ is contained in the current graph, neither rule can delete a vertex or edge belonging to $J$, because its degree or common-neighbor count in the current graph is at least the corresponding count in $J$.
Induction shows that every such $J$ is contained in the terminal graph.
The terminal graph itself satisfies both conditions, so it is the unique largest one, regardless of deletion order.
\end{proof}

The vertex rule alone produces the $(q-1)$-core.
With both rules, the result is the $q$-truss, allowing disconnected components and excluding isolated vertices: every edge lies in at least $q-2$ triangles~\cite{cohen2008trusses}.
Indeed, an endpoint of an edge with $q-2$ common neighbors has degree at least $q-1$.
For $q'\geq q\geq2$, Proposition~\ref{prop:peeling} also gives $H_{q'}\subseteq H_q$, since the terminal graph for $q'$ satisfies the conditions for $q$.

\begin{figure}[t]
\centering
\scalebox{0.6}{
\begin{tikzpicture}[line width=1.8pt]
\coordinate (A) at (2.5,5);
\coordinate (B) at (0,0);
\coordinate (C) at (5,0);
\coordinate (D) at (-2.5,5);
\coordinate (E) at (2.5,-5);
\coordinate (F) at (7.5,5);

\fill[black] (A) circle (0.25);
\fill[black] (B) circle (0.25);
\fill[black] (C) circle (0.25);
\fill[black] (D) circle (0.25);
\fill[black] (E) circle (0.25);
\fill[black] (F) circle (0.25);

\draw[line width=2.4pt] (A) -- (B);
\draw[line width=2.4pt] (B) -- (C);
\draw[line width=2.4pt] (C) -- (A);
\draw (A) -- (D);
\draw (B) -- (D);
\draw (B) -- (E);
\draw (C) -- (E);
\draw (A) -- (F);
\draw (C) -- (F);

\node at (2.5,5.7) {\Large 1 ($0.\overline{96}$)};
\node at (7.5,5.7) {\Large 6 ($0.\overline{60}$)};
\node at (-2.5,5.7) {\Large 4 ($0.\overline{60}$)};
\node at (2.5,-5.7) {\Large 5 ($0.\overline{12}$)};
\node at (-1.4,0) {\Large 2 ($0.\overline{84}$)};
\node at (6.4,0) {\Large 3 ($0.\overline{84}$)};
\end{tikzpicture}
}
\caption{The graph reduction example. Vertex labels give the vertex index and, in parentheses, its value in an optimal solution of the LP relaxation of~\eqref{objN}--\eqref{ConN5} for $q=4$.}
\label{fig:E1}
\end{figure}
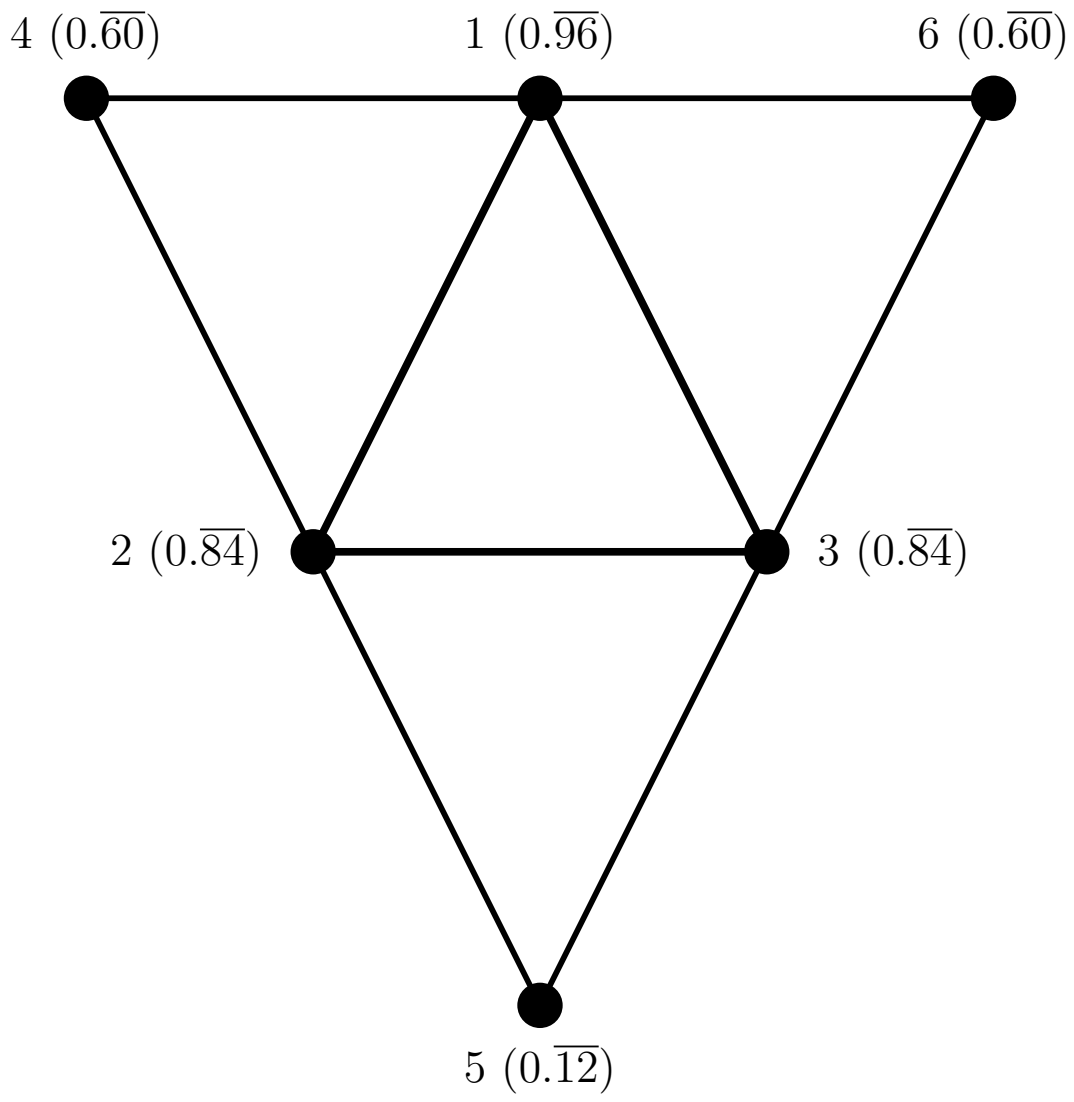

Figure~\ref{fig:E1} illustrates the effect of vertex deletion.
The graph has clique number $3$, while the displayed fractional solution has value $4$ in the LP relaxation of~\eqref{objN}--\eqref{ConN5} at $q=4$.
It is LP-optimal because~\eqref{ConN2} gives $\sum_i y_i\leq4$.
The vertices $4$, $5$, and $6$ have degree $2$ and are deleted first.
The remaining triangle then has degree $2$ at each vertex and is also deleted, certifying that no four-vertex clique exists.

The distinction between the original and reduced auxiliary graphs matters for cut validity.
Because edges as well as vertices are deleted, $H_q$ need not be an induced subgraph of $\bar G$; in particular, $\overline{H_q}$ need not equal $G^k[V_q]$.
A cut valid for every clique of $H_q$ is valid for the fixed-$q$ problem by Proposition~\ref{prop:peeling}, but need not be valid for every clique of the original auxiliary graph.
For example, if two triangles are joined by a bridge $uv$, the reduction for $q=3$ deletes the bridge and keeps both triangles.
The inequality $y_u+y_v\leq1$ is then valid for the reduced graph, although $\{u,v\}$ is a clique of the original graph.

\section{Combinatorial cuts}\label{sec:vi}

Throughout this section, $H=(U,F)$ denotes the current auxiliary graph.
For $W\subseteq U$, write $y(W)=\sum_{i\in W}y_i$.
A clique intersects $W$ in at most $\omega(H[W])=\alpha_1(\overline H[W])$ vertices.
Thus any certified upper bound $r(W)\geq\omega(H[W])$ gives a valid inequality
\begin{equation}\label{eq:rank-certificate}
y(W)\leq r(W).
\end{equation}
When the right-hand side is exact, this is a rank inequality for the stable set polytope of $\overline H$; using an upper bound gives a valid relaxation of that inequality~\cite{chvatal1975rank,rossi2001branchcut,rebennack2011branchcut}.
For $H=\bar G$, the complement is $G^k$; after reduction, validity for the target problem follows from Proposition~\ref{prop:peeling}.
The families below differ in the combinatorial certificate for $r(W)$ and in whether coefficients are lifted.
Section~\ref{sec:sb} replaces the combinatorial certificate with a spectral one.

\begin{proposition}[Scaling a valid inequality by $z$]\label{prop:z-scaling}
Suppose that $c\in\mathbb R^U$ and $r\geq0$ satisfy $\sum_{i\in C}c_i\leq r$ for every $q$-clique $C$ of $H$.
Then
\begin{equation}\label{eq:z-scaled-cut}
\sum_{i\in U}c_i y_i\leq rz
\end{equation}
is valid for~\eqref{objN}--\eqref{ConN5} and dominates $\sum_i c_i y_i\leq r$ on $0\leq z\leq1$.
In particular, $y(W)\leq r(W)z$ is valid whenever $r(W)\geq\omega(H[W])$.
\end{proposition}
\begin{proof}
By Proposition~\ref{prop:fixed-size-exactness}, a feasible solution is either $(0,0)$, for which~\eqref{eq:z-scaled-cut} is an equality, or $(\mathbf1_C,1)$ for a $q$-clique $C$, for which it is the assumed inequality.
Since $r\geq0$ and $z\leq1$, its right-hand side is at most $r$.
This also covers $r=0$ and an empty working graph.
\end{proof}

We present both the unscaled inequalities used in the computational study and the strengthened forms supplied by Proposition~\ref{prop:z-scaling}.
The reported experiments do not evaluate the $z$-scaled variants.
For separation, let $(\bar y,\bar z)$ be the current LP solution.
If $\bar z=0$, then~\eqref{ConN2} and nonnegativity force $\bar y=0$, so none of the inequalities below is violated.

\subsection{Anticlique inequalities}\label{sec:anticlique}

Two nonadjacent vertices cannot belong to the same clique.
This gives the pairwise inequalities
\begin{equation}\label{pi}
y_i+y_j\leq1\qquad
\bigl(i,j\in U,\ i\neq j,\ \{i,j\}\notin F\bigr),
\end{equation}
which underlie the standard edge formulation on the complement graph~\cite{wu2015review}.
For example, the LP solution in Figure~\ref{fig:E1} has
$y_3+y_4=0.\overline{84}+0.\overline{60}>1$, although vertices $3$ and $4$ are nonadjacent.
More generally, an \emph{anticlique} is an independent set $S$ of $H$.
It gives
\begin{equation}\label{ac}
y(S)\leq1\qquad(S\subseteq U\text{ independent in }H),
\end{equation}
the classical clique inequality for the stable set polytope of $\overline H$~\cite{padberg1973facial,grotschel1988geometric}.
The corresponding strengthened forms are
\begin{align}
y_i+y_j&\leq z
&&\bigl(i,j\in U,\ i\neq j,\ \{i,j\}\notin F\bigr),\label{pi-scaled}\\
y(S)&\leq z
&&(S\subseteq U\text{ independent in }H).\label{ac-scaled}
\end{align}
For singleton sets,~\eqref{ac-scaled} is already implied by~\eqref{eq:y-bounded-by-z}.

Exact separation requires a maximum-weight independent set computation in $H$, with weights $\bar y_i$, and is \textsf{NP}-hard in general.
We therefore use a greedy procedure to generate candidate sets.
Sort the vertices by nonincreasing $\bar y_i$.
For a starting vertex $i$ with $\bar y_i>0$, choose the first subsequent, unmarked vertex $j$ that is nonadjacent to $i$; skip $i$ if no such vertex exists.
Starting with $S=\{i,j\}$, scan the remaining vertices in this order and add a vertex whenever it is nonadjacent to every vertex already in $S$.
The experimental implementation adds~\eqref{ac} when $\bar y(S)>1$ and marks the vertices of the accepted set to suppress duplicate candidates within the same separation call.
The same candidate-generation procedure separates~\eqref{ac-scaled} by using the threshold $\bar y(S)>\bar z$.

Anticlique inequalities were separated only within the branch-and-bound tree in the reported experiments.
The implementation also contains an optional preprocessing routine that seeks large anticliques by solving maximum clique problems on $G^k$ under a time limit; this routine was disabled in those experiments.

\subsection{Lifted rank inequalities}\label{sec:cvi}

Fix an integer $a$ with $1\leq a\leq q-1$.
We construct a set $S\subseteq U$ with an insertion order in which each vertex has at most $a$ neighbors among its predecessors.
A greedy coloring in this order uses at most $a+1$ colors, so
$\omega(H[S])\leq a+1$ and $y(S)\leq a+1$ is valid.
This unlifted inequality is implied by the sum of the anticlique inequalities for the color classes; the sum may be stronger if fewer than $a+1$ colors suffice.

Candidate sets are generated from the same ordering by nonincreasing LP value.
For a starting vertex $i$, choose the first subsequent vertex $j$ adjacent to $i$ and initialise $S=\{i,j\}$; if no such $j$ exists, skip $i$.
Scan the remaining vertices and add $u$ whenever $|N_H(u)\cap S|\leq a$.
The initial edge is compatible with this insertion rule because $a\geq1$.

\begin{proposition}[An independent-set lifting]\label{prop:rank-lifting}
Let $S\subseteq U$ satisfy $\omega(H[S])\leq a+1$, where $a\geq1$, and put $b_i=|N_H(i)\cap S|$ for $i\in S$.
Let $L\subseteq\{i\in S:b_i<a\}$ be independent in $H$ and define
\[
c_i=\begin{cases}
a+1-b_i,&i\in L,\\
1,&i\in S\setminus L.
\end{cases}
\]
Then the lifted rank inequality
\begin{equation}\label{cvi}
\sum_{i\in S}c_i y_i\leq a+1
\end{equation}
is valid for every clique of $H$.
Consequently,
\begin{equation}\label{cvi-scaled}
\sum_{i\in S}c_i y_i\leq(a+1)z
\end{equation}
is valid for the fixed-size formulation.
\end{proposition}
\begin{proof}
A clique $C$ contains at most one vertex of $L$.
If $C\cap L=\varnothing$, its left-hand side in~\eqref{cvi} is $|C\cap S|\leq a+1$.
Otherwise, write $C\cap L=\{i\}$.
Every other vertex of $C\cap S$ is a neighbor of $i$ in $S$ and has coefficient $1$.
There are at most $b_i$ such vertices, so the left-hand side is at most
$(a+1-b_i)+b_i=a+1$.
The scaled inequality follows from Proposition~\ref{prop:z-scaling}.
\end{proof}

To choose $L$, repeatedly select an eligible vertex $i$ that is nonadjacent to all previously selected vertices and maximises $(a-b_i)\bar y_i$.
Stop when no eligible vertex remains.
This quantity is the increase in the left-hand side obtained by replacing coefficient $1$ with $a+1-b_i$.
The construction permits $b_i=0$; such a vertex receives coefficient $a+1$, as required when it is the only selected vertex from $S$.
More generally, any certified bound $u_i\geq\omega(H[N_H(i)\cap S])$ can replace $b_i$ in the proof; the resulting coefficient improves on $1$ only when $u_i<a$.
In particular, the coefficient $a+1-\omega(H[N_H(i)\cap S])$ is valid and is at least as large as the degree-based coefficient, although computing it exactly requires a further clique computation.

The experimental separation routine repeats the construction for starting vertices with positive LP value and increases $a$ from $1$, subject to $a+1\leq q$, until the routine finds no further violated inequality at the current value of $a$.
This is a heuristic stopping rule and does not certify that larger values of $a$ yield no violated cuts.
The same candidates can be tested against~\eqref{cvi-scaled} by comparing their weighted LP value with $(a+1)\bar z$.
For $q=1$, this family is unnecessary.

\subsection{Shortest-path inequalities}\label{sec:sp}

\begin{figure}[t]
\centering
\scalebox{0.6}{
\begin{tikzpicture}[line width=1.8pt]
\coordinate (A) at (0,-1);
\coordinate (B) at (0,2);
\coordinate (C) at (2,4.75);
\coordinate (D) at (5,6);
\coordinate (E) at (8,6);

\fill[black] (A) circle (0.25);
\fill[black] (B) circle (0.25);
\fill[black] (C) circle (0.25);
\fill[black] (D) circle (0.25);
\fill[black] (E) circle (0.25);

\draw[line width=2.4pt] (A) -- (B) -- (C) -- (D) -- (E);

\draw[dotted] (A) -- (C);
\draw[dotted] (A) -- (D);
\draw[dotted] (A) -- (E);
\draw[dotted] (B) -- (D);
\draw[dotted] (B) -- (E);
\draw[dotted] (C) -- (E);

\node at (-1,-1) {\Large 1};
\node at (-1,2) {\Large 2};
\node at (1.25,5.5) {\Large 3};
\node at (5,7) {\Large 4};
\node at (8,7) {\Large 5};
\end{tikzpicture}
}
\caption{A shortest path in the auxiliary graph. Solid lines are edges; dotted lines indicate the absent chords. A clique contains at most two of the displayed vertices.}
\label{fig:SP}
\end{figure}
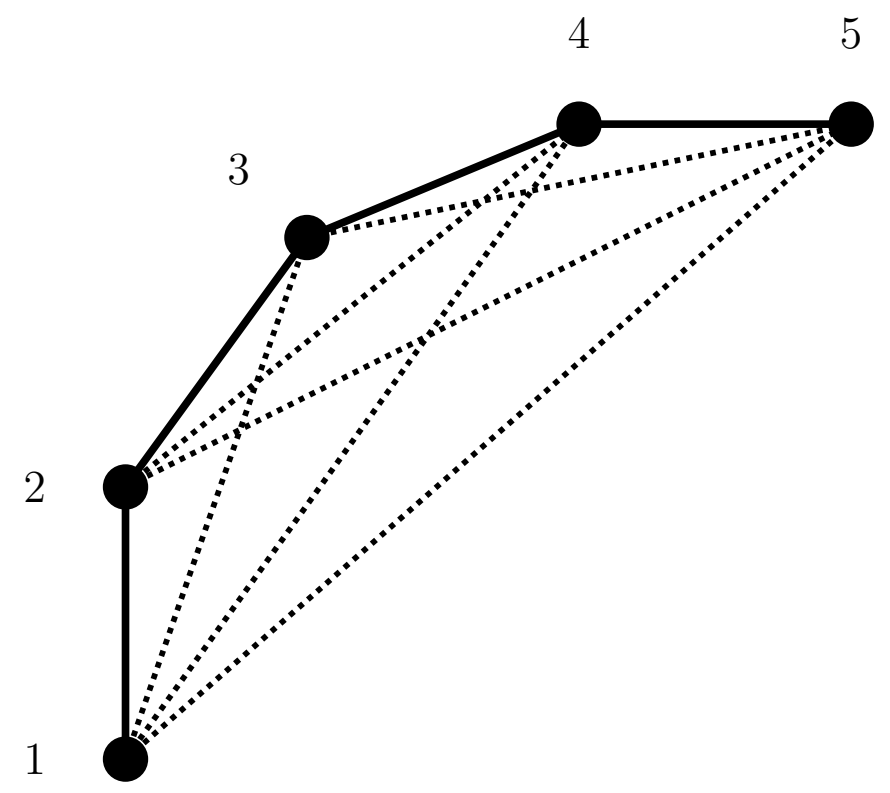

Shortest paths provide anticliques without an LP solution and therefore give inexpensive preprocessing cuts.
Let $P=SP_{\bar G}(i,j)$ be any shortest path in the original auxiliary graph, with
$2\leq\dist_{\bar G}(i,j)<\infty$.
A shortest path in an unweighted graph is induced: an edge joining nonconsecutive path vertices would produce a shorter path.
Consequently, a clique meets $P$ in at most two vertices, and
\begin{equation}\label{SPcut}
\sum_{\ell\in P}y_\ell\leq2
\end{equation}
is valid.
When the working graph has been reduced, variables for deleted vertices are interpreted as zero.
The inequality remains valid because the working graph is a subgraph of $\bar G$.
Figure~\ref{fig:SP} gives the five-vertex example.

Writing the vertices of $P$ in path order, let $S_0$ and $S_1$ be its two parity classes.
Both are independent in $\bar G$, and remain so after any graph reduction.
Thus~\eqref{SPcut} is the sum of $y(S_0)\leq1$ and $y(S_1)\leq1$.
In particular, for a path of length two it follows from the pairwise inequality on the endpoints and the upper bound on the middle vertex.
The scaled path inequality is
\begin{equation}\label{SPcut-scaled}
\sum_{\ell\in P}y_\ell\leq2z,
\end{equation}
and it is dominated by the two separate inequalities
\begin{equation}\label{eq:path-parity-cuts}
y(S_0)\leq z,\qquad y(S_1)\leq z.
\end{equation}
Retaining the parity inequalities therefore gives a formulation at least as strong as retaining only their sum.
It uses twice as many rows, with the same total number of coefficients on the $y$-variables.
The reported implementation generates the unscaled inequalities~\eqref{SPcut} once, before the first subproblem, using one shortest path per eligible vertex pair and only for instances with $|V|\leq600$.
Neither the parity refinement nor the $z$-scaled forms are included in the reported computational results.

\section{Spectral bounds and cuts}\label{sec:spectral}

We use spectral information in two ways: to bound the target clique size globally and to certify the right-hand sides of rank inequalities on subsets of the working graph. We first recall the bounds and specify the numerical conditions needed for valid certificates. We then describe spectral rank cuts and a further inequality for the fixed-size problem. The strengthened inequalities presented below are distinguished from the forms used in the reported experiments.

\subsection{Spectral bounds}\label{sec:spectralboundskMIS}

For a simple graph $J$ on $n$ vertices, let $A$ be its adjacency matrix, with eigenvalues $\lambda_1\geq\cdots\geq\lambda_n$. Let $D$ be the diagonal matrix of vertex degrees and let $\delta$ be the minimum degree. The Laplacian $L=D-A$ has eigenvalues $0=\mu_1\leq\cdots\leq\mu_n$. When $\delta\geq1$, the normalized Laplacian $\mathcal L=I-D^{-1/2}AD^{-1/2}$ is defined and has eigenvalues $0=\sigma_1\leq\cdots\leq\sigma_n$.

\begin{theorem}[Inertia bound~\cite{cvetkovic1971graphs}]\label{thm:inertia_bound}
For a graph $J$ on $n$ vertices,
\[
\alpha_1(J)\leq
\min\bigl\{|\{i:\lambda_i\geq0\}|,\ |\{i:\lambda_i\leq0\}|\bigr\}.
\]
\end{theorem}

\begin{theorem}[Haemers bound~\cite{haemers1995interlacing}]\label{thm:haemers_bound}
For a graph $J$ on $n$ vertices with at least one edge and minimum degree $\delta$,
\[
\alpha_1(J)\leq n\frac{-\lambda_1\lambda_n}{\delta^2-\lambda_1\lambda_n}.
\]
\end{theorem}

\begin{theorem}[Godsil--Newman bound~\cite{godsilnewman2008independent}]\label{thm:godsil_newman_bound}
For a graph $J$ on $n$ vertices with at least one edge and minimum degree $\delta$, the largest Laplacian eigenvalue satisfies $\mu_n>0$, and
\[
\alpha_1(J)\leq n\frac{\mu_n-\delta}{\mu_n}.
\]
\end{theorem}

\begin{theorem}[Harant--Richter bound~\cite{harant_richter2015independent}]\label{thm:harant_richter_bound}
For a graph $J$ on $n$ vertices with $m$ edges and minimum degree $\delta\geq1$,
\[
\alpha_1(J)\leq \frac{2m}{\delta}\frac{\sigma_n-1}{\sigma_n}.
\]
\end{theorem}

The inertia bound uses the signs of the adjacency eigenvalues. The other three bounds use extremal eigenvalues and all reduce to Hoffman's bound $n(-\lambda_n)/(d-\lambda_n)$ when $J$ is $d$-regular with $d>0$. For an edgeless graph, $\alpha_1(J)=n$ directly. More generally, if $I$ is the set of isolated vertices, then
\[
\alpha_1(J)=|I|+\alpha_1(J-I).
\]
We may therefore remove isolated vertices before applying the bounds and add their number afterwards; all parameters of a bound on $J-I$ are evaluated on that graph.

\paragraph{A global bound for \texorpdfstring{$k$-MIS}{k-MIS}}
The identity $\alpha_k(G)=\alpha_1(G^k)=\omega(\bar G)$ allows us to apply Theorems~\ref{thm:inertia_bound}--\ref{thm:harant_richter_bound} to the original power graph $G^k$, before target-dependent reduction of the auxiliary graph. Let $I$ be the isolated vertices of $G^k$ and put $J=G^k-I$. If $J$ is empty, set $q_{\max}=|V|$. Otherwise, let $B_1(J),\ldots,B_4(J)$ be the four bounds above and define, in exact arithmetic,
\begin{equation}\label{eq:spectral-qmax}
q_{\max}=\min\left\{|V|,\ |I|+\left\lfloor\min_{1\leq j\leq4}B_j(J)\right\rfloor\right\}.
\end{equation}
This is an integer upper bound on $\alpha_k(G)$. Bounds based directly on the spectrum of $G$, including polynomial inertia-type and ratio-type bounds~\cite{abiad2016spectral,abiad2019kindependence,ACFNS2022,abiad2023eigenvalue}, can be incorporated by taking the minimum with any additional certified upper bound. Their computation and their effect on the algorithm are separate from the four-bound construction in~\eqref{eq:spectral-qmax}.

\paragraph{Numerical certification}
The validity of a spectral cut or stopping bound depends on the direction of its numerical error. The following conditions provide sufficient certificates; they are theoretical requirements and do not assert that the reported runs were validated by interval computations. Here all spectra and graph parameters refer to the graph being bounded, so $n=|V(J)|$ when isolated vertices have been removed.

Suppose intervals $[a_i,b_i]$ enclose the adjacency eigenvalues, counted with multiplicity. Then
\begin{equation}\label{eq:certified-inertia}
B_{\mathrm{in}}^{+}(J):=
\min\bigl\{|\{i:b_i\geq0\}|,\ |\{i:a_i\leq0\}|\bigr\}
\end{equation}
is a valid inertia upper bound. Indeed, every nonnegative eigenvalue is counted in the first set and every nonpositive eigenvalue in the second; an interval containing zero is conservatively counted in both. In particular, the usual sign-tolerance rule is valid when the tolerance is a certified bound on every eigenvalue error.

For the ratio bounds, suppose positive numbers $u,v,M,S$ satisfy
\[
\lambda_1\leq u,\qquad -\lambda_n\leq v,\qquad
\mu_n\leq M,\qquad \sigma_n\leq S.
\]
Under the respective theorem hypotheses, conservative upper bounds are
\begin{equation}\label{eq:certified-ratio-bounds}
\frac{nuv}{\delta^2+uv},\qquad
n\left(1-\frac{\delta}{M}\right),\qquad
\frac{2m}{\delta}\left(1-\frac1S\right).
\end{equation}
For the first expression, $-\lambda_1\lambda_n\leq uv$ and $t\mapsto nt/(\delta^2+t)$ is nondecreasing for $t>0$. The other two expressions are nondecreasing in the positive arguments $M$ and $S$. Evaluating these expressions with outward rounding and using the upper endpoint of each resulting enclosure therefore preserves validity. Integer rounding is applied only after this upper endpoint has been obtained. An uncertified candidate may be omitted without affecting the validity of the remaining bounds.

A fixed additive guard on the final bound is not a substitute for this calculation. For example, the convention of treating eigenvalues within $10^{-4}$ of zero as zero and rounding a computed ratio bound $\widehat\beta$ to $\lfloor\widehat\beta+0.01\rfloor$ does not by itself certify the latter bound. A matching on 1000 vertices has $\alpha_1=500$ and $\sigma_n=2$. Substituting $\widehat\sigma_n=1.99991$ in Theorem~\ref{thm:harant_richter_bound} gives $\widehat\beta\approx499.977499$ and hence $\lfloor\widehat\beta+0.01\rfloor=499$, despite an eigenvalue error smaller than $10^{-4}$. Certification is required for every spectral cut used in a proof, as well as for a global bound used to terminate the search; validating only the final stopping bound does not validate earlier cuts.

\subsection{Spectral bound inequalities}\label{sec:sb}

Fix a target size $q$ and let $H=(U,F)$ be the current working auxiliary graph. After peeling, $H=H_q$ in the notation of Section~\ref{sec:rgs}. Let $(\bar y,\bar z)$ be the current LP solution, and order $U$ as $\pi_1,\ldots,\pi_{|U|}$ so that $\bar y_{\pi_1}\geq\cdots\geq\bar y_{\pi_{|U|}}$. For $2\leq\ell\leq|U|$, put
\[
W_\ell=\{\pi_1,\ldots,\pi_\ell\},\qquad
\widetilde G_\ell=\overline{H[W_\ell]},
\]
and let $A_\ell$ be the adjacency matrix of $\widetilde G_\ell$. The complement here is taken within $W_\ell$. Before edge peeling, $\widetilde G_\ell=G^k[W_\ell]$; after peeling this equality need not hold.

The intersection of a clique of $H$ with $W_\ell$ has size at most $\alpha_1(\widetilde G_\ell)$. Thus any upper bound on this independence number gives a rank inequality for the stable set polytope of $\overline H$~\cite{chvatal1975rank,rossi2001branchcut,rebennack2011branchcut}. Using the inertia bound, define
\begin{equation}\label{sb}
\beta_\ell:=\min\bigl\{|\{i:\lambda_i(A_\ell)\geq0\}|,\
|\{i:\lambda_i(A_\ell)\leq0\}|\bigr\}
\geq\alpha_1(\widetilde G_\ell).
\end{equation}
The resulting spectral rank cut is
\begin{equation}\label{sb-cut}
\sum_{i=1}^{\ell}y_{\pi_i}\leq\beta_\ell.
\end{equation}
When exact inertia is unavailable, $\beta_\ell$ in a cut is replaced by a certified integer upper bound such as~\eqref{eq:certified-inertia}.

We strengthen this spectral rank cut for the fixed-size formulation to
\begin{equation}\label{sb-cut-scaled}
\sum_{i=1}^{\ell}y_{\pi_i}\leq\beta_\ell z.
\end{equation}
At an integer feasible solution, $z=0$ forces $y=0$, while $z=1$ selects a $q$-clique; these are the only possibilities even when $z$ is declared continuous, as shown in Section~\ref{sec:MCP}. Hence~\eqref{sb-cut-scaled} is valid. Since $0\leq z\leq1$, it dominates~\eqref{sb-cut}. The reported spectral rank cuts use the unscaled form~\eqref{sb-cut}; the computational effect of~\eqref{sb-cut-scaled} has not been measured in those experiments.

The inertia certificate can improve on a coloring certificate. For example, when $\widetilde G_\ell=C_5$, its inertia bound is 2 whereas $\chi(\overline{C_5})=3$. Such an improvement is not implied by symmetry alone. The Lov\'asz theta number of $\widetilde G_\ell$ provides another upper bound on its independence number~\cite{lovasz1979shannon}. A certified theta bound could be combined with inertia by taking their minimum, but would require a semidefinite optimization step; neither certificate uniformly dominates the other.

\paragraph{Separation and pruning}
Write $s_\ell=\sum_{i=1}^{\ell}\bar y_{\pi_i}$. The unscaled cut is violated when $s_\ell>\beta_\ell$. A dense eigendecomposition of $A_\ell$ requires $O(\ell^3)$ arithmetic operations, so checking every prefix independently costs $O(|U|^4)$ operations per separation call. The separator therefore uses filters to reduce the number of decompositions.

One filter skips a prefix when $\bar y_{\pi_{\ell-1}}+\bar y_{\pi_\ell}\leq1$. This is a heuristic and does not certify nonviolation. For example, let $H$ consist of a $C_5$ and an additional vertex adjacent to all five cycle vertices, and take $q=3$. This graph survives the reduction rules. Assigning LP value $1/2$ to each cycle vertex, zero to the additional vertex, and $\bar z=5/6$ gives a feasible point of the base LP. The cycle vertices form a prefix with $s_5=5/2>\beta_5=2$, since $\overline{C_5}\cong C_5$, but the test skips it.

An exact pruning rule follows from interlacing. Since $A_\ell$ is a principal submatrix of $A_{\ell+1}$, Cauchy interlacing implies that the latter has at least as many nonnegative eigenvalues and at least as many nonpositive eigenvalues as the former. Consequently,
\[
\beta_{\ell+1}\geq\beta_\ell.
\]
Once $\beta_\ell$ is known exactly, a later prefix $j>\ell$ cannot violate~\eqref{sb-cut} while $s_j\leq\beta_\ell$. Computations may therefore be postponed until the first later prefix whose mass exceeds this value. Independently computed conservative upper bounds on the inertia values need not be monotone. If exact values are unavailable, this safe pruning test instead requires a certified lower bound on $\beta_\ell$; using an arbitrary upper bound gives only a heuristic filter. Skipping a cut heuristically affects separation strength, not the validity of cuts actually added.

Two further comparisons avoid retaining unnecessary cuts. If $\ell<|U|$ and $\beta_\ell=\beta_{\ell+1}$, the cut for $\ell+1$ dominates the cut for $\ell$ by nonnegativity. Also, if $\ell_1<\ell$ and
\[
s_\ell-s_{\ell_1}\leq\beta_\ell-\beta_{\ell_1},
\]
then $s_{\ell_1}-\beta_{\ell_1}\geq s_\ell-\beta_\ell$: the earlier cut has at least as much violation at the current LP point. The latter comparison is not a polyhedral dominance statement. A cut is discarded in favor of another only when the preferred inequality is retained or already present. With numerical certificates, these two comparisons use the actual right-hand sides of the certified cuts.

For the strengthened family~\eqref{sb-cut-scaled}, violation is tested against $\beta_\ell\bar z$. The safe prefix test becomes $s_j\leq\beta_\ell\bar z$, and the comparison of current violations becomes
\[
s_\ell-s_{\ell_1}\leq
(\beta_\ell-\beta_{\ell_1})\bar z.
\]
The same certification qualifications apply. In particular, the unscaled pruning threshold cannot be used as a nonviolation certificate for the scaled family when $\bar z<1$.

\paragraph{Extension and lifting}
Suppose a set $W\subseteq U$ has a certified integer bound $r\geq\omega(H[W])$. A vertex outside $W$ may be added if it has at most $r-1$ neighbors in the current set. Every clique containing the added vertex then has size at most $r$, while cliques not containing it already satisfy the bound. Applying this argument after each insertion preserves $\omega(H[W])\leq r$. Candidates may be considered in nonincreasing LP order.

The coefficient lifting of Section~\ref{sec:cvi} also applies with right-hand side $r$. Specifically, for a pairwise nonadjacent set of lifted vertices in $W$, a vertex $i$ with $d_i=|N_H(i)\cap W|<r-1$ may receive coefficient $r-d_i$, while all unlifted coefficients remain 1. A clique contains at most one lifted vertex; if it contains $i$, the remaining contribution is at most $d_i$, so the total is at most $r$. If it contains none, its contribution is its cardinality, also at most $r$. The resulting inequality may again be strengthened by replacing its right-hand side $r$ with $rz$.

\subsection{Non-selected vertex inequalities}\label{sec:nonSelected}

This family exploits a different implication of the fixed target size: a $q$-clique containing sufficiently many vertices of a chosen set may have to use vertices outside the support of the current LP solution. For any $T\subseteq U$, write $y(T)=\sum_{i\in T}y_i$.

For each integer $1\leq a\leq q-1$, scan the vertices with $0<\bar y_i<1$ in nonincreasing LP order and let $S$ be the first prefix with $\bar y(S)>a$. If their total mass is at most $a$, skip this value of $a$. A vertex with $\bar y_i=1$ need not be included in the candidate set: the base neighborhood constraint requires neighbor mass at least $q-1$, while the cardinality equation and $\bar z\leq1$ leave at most $q-1$ mass outside $i$. Thus $\bar z=1$ and every other vertex with positive LP value is adjacent to $i$.

Obtain a certified integer bound $b\geq\omega(H[S])$, for example from the inertia of $\overline{H[S]}$. Since $S$ is nonempty, $b\geq1$. Continue only when $b\leq a$, and define
\begin{align*}
P&=\{i\in U\setminus S:\bar y_i>0,\ |N_H(i)\cap S|\geq b\},\\
R&=\{i\in U\setminus S:\bar y_i=0,\ |N_H(i)\cap S|\geq b\}.
\end{align*}
Let $c\geq\omega(H[P])$ be a certified integer bound computed on the entire set $P$, with $c=0$ for $P=\varnothing$ and $c=1$ for $|P|=1$. If $q\leq b+c$, no cut of the following form is generated.

\begin{proposition}\label{prop:nonselected-validity}
Under the preceding assumptions, if $q>b+c$, the inequality
\begin{equation}\label{nsv-scaled}
\bigl(y(S)-(b-1)z\bigr)(q-b-c)\leq y(R)
\end{equation}
is valid for the integer feasible solutions of the fixed-size formulation.
\end{proposition}

\begin{proof}
If $z=0$, the cardinality equation and nonnegativity force $y=0$, so both sides are zero. Otherwise $z=1$, and the selected vertices form a $q$-clique $C$ in $H$. The bound on $H[S]$ gives $|C\cap S|\leq b$. If $|C\cap S|\leq b-1$, the left-hand side is nonpositive and the inequality follows from $y(R)\geq0$.

If $|C\cap S|=b$, every vertex of $C\setminus S$ is adjacent to these $b$ vertices and therefore lies in $P\cup R$. At most $c$ selected vertices belong to $P$. Hence $|C\cap R|\geq q-b-c$, which is precisely~\eqref{nsv-scaled} in this case.
\end{proof}

The corresponding unscaled inequality is
\begin{equation}\label{nsv}
\bigl(y(S)-b+1\bigr)(q-b-c)\leq y(R).
\end{equation}
It follows from~\eqref{nsv-scaled} because $b\geq1$ and $z\leq1$; the difference between their left-hand sides is $(b-1)(1-z)(q-b-c)\geq0$. We thus obtain a scaled inequality that is at least as strong throughout the LP domain. Under the separation condition $b\leq a<\bar y(S)$, the ordinary rank cut $y(S)\leq b$ is already violated and can also be added.

The requirement that $c$ bound the clique number of the full set $P$ is essential. A bound $c_0$ computed on a restricted set $P_0\subseteq P$ can be extended safely by taking $c=c_0+|P\setminus P_0|$, since every clique uses at most $c_0$ vertices of $P_0$ and at most $|P\setminus P_0|$ others. Using $c_0$ alone is not justified.

The reported experiments use the earlier inequality
\begin{equation}\label{nsv-legacy}
\bigl(y(S)-b\bigr)(q-b-c)\leq y(R),\qquad q>b+c,
\end{equation}
rather than~\eqref{nsv} or~\eqref{nsv-scaled}. In the presence of the rank cut $y(S)\leq b$,~\eqref{nsv-legacy} is redundant: its left-hand side is nonpositive and its right-hand side is nonnegative. The corrected and scaled forms above are analytical improvements whose computational effect is not established by those experiments.

\section{Solution procedure}\label{sec:solutionProceedure}

We solve a sequence of fixed-size clique problems on the auxiliary graph $\bar G=\overline{G^k}$. Starting with $q=1$, it searches for a clique of size $q$, increases the target after a successful search, and stops when it establishes that the next target is unattainable or reaches a valid upper bound $q_{\max}$. This order favors successful feasibility searches before the first unsuccessful target. Its computational motivation is the concentration of effort near that target observed in Section~\ref{sec:effectQ}; it is not a worst-case advantage over binary search.

For each target $q\geq2$, the reduction of Section~\ref{sec:rgs} produces $H_q=(V_q,E_q)$. The fixed-size model is then constructed on $H_q$. A feasible solution selecting $q$ vertices attains the model's upper objective bound $q$, so the search advances as soon as such a solution is found. If graph reduction rules out a $q$-clique, or the solver establishes optimal objective value zero, the previous successful target gives $\alpha_k(G)=q-1$. Alternatively, finding a clique of size $q_{\max}$ gives $\alpha_k(G)=q_{\max}$. The fixed-size model itself remains feasible at the zero solution throughout.

\paragraph{Separation and cut reuse}
Shortest-path inequalities are generated during preprocessing in the unreduced auxiliary graph and are retained when the target increases. At separation callbacks, vertices are sorted by nonincreasing current LP value, and the procedures of Sections~\ref{sec:vi} and~\ref{sec:spectral} search for violated cuts. Completeness of these separation procedures is unnecessary: the base integer model is already exact, and valid cuts improve its relaxation without removing feasible cliques.

Rank inequalities and their liftings can also be carried from one target to the next. The reductions are monotone in the target: for $q'>q$, the graph $H_{q'}$ is a subgraph of $H_q$. An inequality valid for every clique of $H_q$ therefore remains valid on $H_{q'}$, with variables for deleted vertices fixed to zero. The same argument justifies reusing the strengthened inequalities with a right-hand side multiplied by $z$: at $z=0$ all selected-vertex variables vanish, and at $z=1$ the underlying clique inequality applies. The non-selected-vertex inequalities depend explicitly on $q$ and must instead be regenerated for the new target.

\begin{theorem}[Correctness of the sequential framework]\label{thm:framework-correctness}
Let $G$ be a nonempty graph and let $q_{\max}\in\{1,\ldots,|V|\}$ satisfy $\alpha_k(G)\leq q_{\max}$. Suppose that, at each target $q$, the reduction preserves every $q$-clique of $\bar G$, and every added inequality is valid for the empty solution and for the incidence vector of each such surviving clique. If each fixed-size subproblem is solved to optimality, the sequential framework returns $\alpha_k(G)$. Complete separation of any cut family is not required. If a subproblem is interrupted, the last clique found supplies a valid lower bound.
\end{theorem}

\begin{proof}
Consider a nonzero feasible solution of the fixed-size model, with binary $y$ and $0\leq z\leq1$. For any selected vertex, \eqref{ConN1} requires at least $q-1$ selected neighbors. Thus at least $q$ vertices are selected, whereas \eqref{ConN2} allows at most $q$. Exactly $q$ vertices are therefore selected, $z=1$, and every selected vertex is adjacent to all the others. Conversely, the incidence vector of a $q$-clique, together with $z=1$, satisfies the base model and, by hypothesis, every added inequality. Reduction preserves these cliques. Hence the fixed-size optimum is $q$ if a $q$-clique exists in $\bar G$, and zero otherwise.

If the first zero optimum occurs at target $q$, a clique of size $q-1$ has already been found. No larger clique exists, since every clique of size at least $q$ contains a $q$-clique. The optimum is therefore $q-1$. If every target through $q_{\max}$ succeeds, the clique of size $q_{\max}$ meets the valid upper bound. In either case the returned value equals $\omega(\bar G)=\alpha_k(G)$. An interrupted search does not invalidate a previously found clique.
\end{proof}

The theorem separates the mathematical framework from its numerical implementation. Spectral bounds and spectral cut coefficients must meet the certification conditions of Section~\ref{sec:spectralboundskMIS} for the theorem to provide an exact certificate. A feasible clique can be checked directly from its vertex set; certifying a successful stopping decision also requires valid upper bounds and cuts throughout the run.

\paragraph{Computational variant evaluated below}
The experiments use the earlier implementation of this framework: $z$ is declared binary, the redundant rows~\eqref{ConN3} are retained, and cuts use their unscaled right-hand sides. The dynamically separated inequalities are \eqref{ac}, \eqref{cvi}, \eqref{sb-cut}, and the weaker non-selected-vertex inequality~\eqref{nsv-legacy}. The corrected inequality~\eqref{nsv}, its strengthened form~\eqref{nsv-scaled}, the strengthened spectral cuts~\eqref{sb-cut-scaled}, and the analogous $z$-scaled combinatorial cuts are not evaluated in those runs. Large-anticlique preprocessing is disabled, and shortest-path preprocessing is used only when $|V|\leq600$.

Within each subproblem, a family is disabled after its separation counter reaches $100$. The counter increases by one at a callback that finds no violated inequality of that family and decreases by ten when a violated inequality is found. Rank cuts from previous targets are retained; target-dependent non-selected-vertex inequalities are not. If the time limit is reached without a $q$-clique or a nonexistence conclusion, the run stops with the lower bound from its last successful target.

\section{Computational results}\label{sec:IPkMIS}

We evaluate the computational variant described above on three benchmark classes for $k\in\{1,2,3,4\}$. We report solution statuses and lower bounds, examine how work is distributed over target sizes, and compare the tested cut configurations. We also compare the method with the discretized integer programming formulation of Martins~\cite{martins2010extended}. The conventions below specify how the recorded numerical outcomes relate to the exact framework.

\subsection{Benchmark instances and experimental conventions}\label{sec:benchmark}

\paragraph{Named graphs}
The first class contains $79$ named graphs from SageMath's graph library~\cite{sagemath}, including the Petersen, Coxeter, and Hoffman--Singleton graphs. Their orders range from $4$ to $600$; the complete list is in Table~\ref{tab:graph-results}.

\paragraph{Coding-theoretic graphs}
The second class contains $32$ graphs from Sloane's challenge on independent sets~\cite{sloane_challenge}. The prefixes \texttt{1dc.*} and \texttt{2dc.*} denote single- and double-deletion-correcting codes; \texttt{1tc.*} denotes single-transposition-correcting codes; \texttt{1et.*} adds end-around transposition; and \texttt{1zc.*} denotes Z-channel codes. The suffix gives the number of vertices. Independence in these graphs represents the corresponding code constraint, and increasing $k$ imposes additional graph-distance separation.

\paragraph{Sum-rank-metric graphs}
The third class contains $25$ sum-rank-metric graphs from Abiad, Khramova, and Ravagnani~\cite{abiad2023eigenvalue}. For a prime power $q$ and tuples $n=(n_1,\ldots,n_t)$ and $m=(m_1,\ldots,m_t)$ of positive integers, the vertices are tuples $(X_1,\ldots,X_t)$ with $X_i\in\mathbb F_q^{n_i\times m_i}$. Two vertices are adjacent when their difference has rank one in exactly one block and is zero in all other blocks. Graph distance equals sum-rank distance, so a $k$-independent set is a code of minimum sum-rank distance at least $k+1$. These graphs have
\[
 |V|=q^{\sum_{i=1}^{t}n_i m_i},\qquad
 |E|=\frac{|V|}{2(q-1)}\sum_{i=1}^{t}(q^{n_i}-1)(q^{m_i}-1).
\]
The field-size symbol $q$ in these instance descriptions is unrelated to the target clique size. Table~\ref{tab:coding-results} lists the parameters explicitly. The Hamming graph is the special case $n_i=m_i=1$ for every $i$; none of the listed instances has all blocks of this form.

\paragraph{Computing environment and limits}
The reported implementation uses C++ with Visual Studio 2017, version 15.9.8, and IBM ILOG CPLEX Optimization Studio 12.9. The runs use one thread on an AMD Ryzen Threadripper 3990X processor at 2.9\,GHz with 256\,GB of RAM. Each fixed-target MIP solve has a three-hour time limit. Total runtime includes successive target searches and preprocessing, and can therefore substantially exceed three hours.

\paragraph{Numerical status}
An asterisk in the tables denotes \emph{solver-reported optimality} for the indicated implementation. For the sequential method, the run terminated with a reported zero optimum at the next target or with a clique attaining its computed upper bound $q_{\max}$; for DF, it records the MIP solver's optimal status. For the configurations using spectral information, rigorous eigenvalue enclosures and a complete audit of the spectral cuts are not available for these runs. The tables therefore distinguish computationally reported optimality from independently certified exact values. The numerical safeguards and certification requirements in Section~\ref{sec:spectralboundskMIS} apply to both global bounds and generated cuts. Values established in the literature are identified separately where relevant.

\paragraph{Table conventions}
Tables~\ref{tab:graph-results}--\ref{tab:coding-results} primarily report the configuration with both cut families, denoted Both. Three entries replace an unresolved Both result by a successful single-family result; these entries carry $^\dagger$ and use that single-family run's time. A dagger thus records an additional observed success, not the runtime of a prescribed portfolio. For an unstarred entry, the reported lower bound is the largest obtained by Both for any $k'\geq k$ on the same graph, using $\alpha_k(G)\geq\alpha_{k'}(G)$. The displayed time remains that of the run at $k$, and need not be the time at which the transferred bound was found. Main-table times are rounded down to whole seconds. Appendix~\ref{app:ablation} reports each configuration's own values without this transfer across $k$, with times rounded to the nearest second.

\subsection{Solution status and lower bounds}\label{sec:effect}

The three main benchmark classes give $544$ instance--$k$ combinations. Both reports optimality for $470$ of them ($86.4\%$). Including the three additional single-family successes increases the union to $473$ ($86.9\%$). Among the $416$ combinations with at most $256$ vertices, these counts are $414$ ($99.5\%$) and $415$ ($99.8\%$), respectively. These percentages use the numerical-status convention of Section~\ref{sec:benchmark}.

For the named graphs, Both reports optimality on $313$ of $316$ combinations. All graphs with at most $100$ vertices finish within seconds for every tested $k$. The three unresolved combinations are Cell120 at $k=1,3,4$; the $k=2$ run returns a starred value of $120$ in $645$ seconds.

For the coding-theoretic graphs, all combinations with at most $128$ vertices are starred. At $256$ vertices, the only unresolved Both run is \texttt{1et.256} at $k=1$; the combinatorial configuration supplies the dagger-marked value $50$. Larger instances distinguish the graph families and target distances. For example, several \texttt{2dc.*} combinations remain tractable as the order increases, whereas many larger instances from the other families retain only lower bounds.

For the sum-rank-metric graphs, all but one combination with at most $256$ vertices are starred; the exception has $t=7$ and $k=2$. At $512$ vertices, all six $k=1$ combinations are starred, but seven of the eighteen combinations with $k\geq2$ are unresolved. All four combinations at $729$ vertices are starred. At $1024$ vertices, $k=1,2$ are starred for Both and $k=4$ is supplied by the spectral-only configuration. At orders at least $2048$, the only starred result comes from the combinatorial-only configuration at $k=1$ on the $2048$-vertex instance. Thus neither the effect of increasing $k$ nor the role of graph order is uniform across the benchmark classes.

The reported values agree with the overlapping exact values for named graphs in~\cite[Table~2]{ACFNS2022} and for sum-rank-metric graphs in~\cite[Tables~1--2]{abiad2023eigenvalue}. In particular, two entries that remain unstarred in our runs are known exactly from the latter source: $\alpha_2=18$ for $t=3$, $q=2$, $n=m=(2,2,1)$, and $\alpha_3=9$ for $t=6$, $q=2$, $n=m=(2,1^5)$. They are marked $^{\S}$ in Table~\ref{tab:coding-results}. Agreement with previously established values provides a consistency check; it does not replace verification of the other numerical conclusions.

\begingroup
	\scriptsize
	\setlength{\tabcolsep}{2pt}
	\renewcommand{\arraystretch}{1.08}

	\begin{center}
		\begin{longtable}{l r r
				!{\,\vrule\,}
				r r !{\,\vrule\,} r r !{\,\vrule\,} r r !{\,\vrule\,} r r}
			\caption{Named graph instances, ordered by $|V|$. Runtimes are in seconds, rounded down. An asterisk indicates solver-reported optimality under the numerical-status convention of Section~\ref{sec:benchmark}; unstarred values are reported lower bounds.}\label{tab:graph-results}\\
			\toprule
			Instance & $|V|$ & $|E|$ &
			$\alpha_1$ & Runtime &
			$\alpha_2$ & Runtime &
			$\alpha_3$ & Runtime &
			$\alpha_4$ & Runtime \\
			\midrule
			\endfirsthead
			
			\toprule
			Instance & $|V|$ & $|E|$ &
			$\alpha_1$ & Runtime &
			$\alpha_2$ & Runtime &
			$\alpha_3$ & Runtime &
			$\alpha_4$ & Runtime \\
			\midrule
			\endhead
			
			\bottomrule
			\endlastfoot
			
			Tetrahedral & 4 & 6 & \alphaopt{1} & 0 & \alphaopt{1} & 0 & \alphaopt{1} & 0 & \alphaopt{1} & 0 \\
			Thomsen & 6 & 9 & \alphaopt{3} & 0 & \alphaopt{1} & 0 & \alphaopt{1} & 0 & \alphaopt{1} & 0 \\
			Octahedral & 6 & 12 & \alphaopt{2} & 0 & \alphaopt{1} & 0 & \alphaopt{1} & 0 & \alphaopt{1} & 0 \\
			MoserSpindle & 7 & 11 & \alphaopt{2} & 0 & \alphaopt{1} & 0 & \alphaopt{1} & 0 & \alphaopt{1} & 0 \\
			Hexahedral & 8 & 12 & \alphaopt{4} & 0 & \alphaopt{2} & 0 & \alphaopt{1} & 0 & \alphaopt{1} & 0 \\
			Wagner & 8 & 12 & \alphaopt{3} & 0 & \alphaopt{1} & 0 & \alphaopt{1} & 0 & \alphaopt{1} & 0 \\
			Petersen & 10 & 15 & \alphaopt{4} & 0 & \alphaopt{1} & 0 & \alphaopt{1} & 0 & \alphaopt{1} & 0 \\
			KrackhardtKite & 10 & 18 & \alphaopt{4} & 0 & \alphaopt{2} & 0 & \alphaopt{2} & 0 & \alphaopt{1} & 0 \\
			Herschel & 11 & 18 & \alphaopt{6} & 0 & \alphaopt{2} & 0 & \alphaopt{2} & 0 & \alphaopt{1} & 0 \\
			Grotzsch & 11 & 20 & \alphaopt{5} & 0 & \alphaopt{1} & 0 & \alphaopt{1} & 0 & \alphaopt{1} & 0 \\
			GoldnerHarary & 11 & 27 & \alphaopt{6} & 0 & \alphaopt{1} & 0 & \alphaopt{1} & 0 & \alphaopt{1} & 0 \\
			Franklin & 12 & 18 & \alphaopt{6} & 0 & \alphaopt{2} & 0 & \alphaopt{1} & 0 & \alphaopt{1} & 0 \\
			BidiakisCube & 12 & 18 & \alphaopt{5} & 0 & \alphaopt{2} & 0 & \alphaopt{1} & 0 & \alphaopt{1} & 0 \\
			Frucht & 12 & 18 & \alphaopt{5} & 0 & \alphaopt{3} & 0 & \alphaopt{2} & 0 & \alphaopt{1} & 0 \\
			Tietze & 12 & 18 & \alphaopt{5} & 0 & \alphaopt{3} & 0 & \alphaopt{1} & 0 & \alphaopt{1} & 0 \\
			TruncatedTetrahedral & 12 & 18 & \alphaopt{4} & 0 & \alphaopt{3} & 0 & \alphaopt{1} & 0 & \alphaopt{1} & 0 \\
			Durer & 12 & 18 & \alphaopt{4} & 0 & \alphaopt{2} & 0 & \alphaopt{2} & 0 & \alphaopt{1} & 0 \\
			Chvatal & 12 & 24 & \alphaopt{4} & 0 & \alphaopt{1} & 0 & \alphaopt{1} & 0 & \alphaopt{1} & 0 \\
			Icosahedral & 12 & 30 & \alphaopt{3} & 0 & \alphaopt{2} & 0 & \alphaopt{1} & 0 & \alphaopt{1} & 0 \\
			Heawood & 14 & 21 & \alphaopt{7} & 0 & \alphaopt{2} & 0 & \alphaopt{1} & 0 & \alphaopt{1} & 0 \\
			Poussin & 15 & 39 & \alphaopt{5} & 0 & \alphaopt{2} & 0 & \alphaopt{1} & 0 & \alphaopt{1} & 0 \\
			Hoffman & 16 & 32 & \alphaopt{8} & 0 & \alphaopt{2} & 0 & \alphaopt{2} & 0 & \alphaopt{1} & 0 \\
			MoebiusKantor & 16 & 24 & \alphaopt{8} & 0 & \alphaopt{4} & 0 & \alphaopt{2} & 0 & \alphaopt{1} & 0 \\
			Clebsch & 16 & 40 & \alphaopt{5} & 0 & \alphaopt{1} & 0 & \alphaopt{1} & 0 & \alphaopt{1} & 0 \\
			Sousselier & 16 & 27 & \alphaopt{6} & 0 & \alphaopt{3} & 0 & \alphaopt{1} & 0 & \alphaopt{1} & 0 \\
			Shrikhande & 16 & 48 & \alphaopt{4} & 0 & \alphaopt{1} & 0 & \alphaopt{1} & 0 & \alphaopt{1} & 0 \\
			Errera & 17 & 45 & \alphaopt{5} & 0 & \alphaopt{2} & 0 & \alphaopt{2} & 0 & \alphaopt{1} & 0 \\
			Pappus & 18 & 27 & \alphaopt{9} & 0 & \alphaopt{3} & 0 & \alphaopt{3} & 0 & \alphaopt{1} & 0 \\
			BlanusaFirstSnark & 18 & 27 & \alphaopt{8} & 0 & \alphaopt{4} & 0 & \alphaopt{2} & 0 & \alphaopt{1} & 0 \\
			BlanusaSecondSnark & 18 & 27 & \alphaopt{7} & 0 & \alphaopt{4} & 0 & \alphaopt{2} & 0 & \alphaopt{1} & 0 \\
			Robertson & 19 & 38 & \alphaopt{7} & 0 & \alphaopt{3} & 0 & \alphaopt{1} & 0 & \alphaopt{1} & 0 \\
			Desargues & 20 & 30 & \alphaopt{10} & 0 & \alphaopt{4} & 0 & \alphaopt{2} & 0 & \alphaopt{2} & 0 \\
			Dodecahedral & 20 & 30 & \alphaopt{8} & 0 & \alphaopt{4} & 0 & \alphaopt{2} & 0 & \alphaopt{2} & 0 \\
			FlowerSnark & 20 & 30 & \alphaopt{9} & 0 & \alphaopt{5} & 0 & \alphaopt{2} & 0 & \alphaopt{1} & 0 \\
			Folkman & 20 & 40 & \alphaopt{10} & 0 & \alphaopt{3} & 0 & \alphaopt{2} & 0 & \alphaopt{1} & 0 \\
			Brinkmann & 21 & 42 & \alphaopt{7} & 0 & \alphaopt{3} & 0 & \alphaopt{1} & 0 & \alphaopt{1} & 0 \\
			Kittell & 23 & 63 & \alphaopt{7} & 0 & \alphaopt{3} & 0 & \alphaopt{2} & 0 & \alphaopt{1} & 0 \\
			Nauru & 24 & 36 & \alphaopt{12} & 0 & \alphaopt{6} & 0 & \alphaopt{4} & 0 & \alphaopt{1} & 0 \\
			McGee & 24 & 36 & \alphaopt{10} & 0 & \alphaopt{5} & 0 & \alphaopt{2} & 0 & \alphaopt{1} & 0 \\
			Markstroem & 24 & 36 & \alphaopt{9} & 0 & \alphaopt{6} & 0 & \alphaopt{3} & 0 & \alphaopt{3} & 0 \\
			Klein7Regular & 24 & 84 & \alphaopt{6} & 0 & \alphaopt{3} & 0 & \alphaopt{1} & 0 & \alphaopt{1} & 0 \\
			F26A & 26 & 39 & \alphaopt{13} & 0 & \alphaopt{6} & 0 & \alphaopt{3} & 0 & \alphaopt{2} & 0 \\
			Holt & 27 & 54 & \alphaopt{10} & 0 & \alphaopt{3} & 0 & \alphaopt{1} & 0 & \alphaopt{1} & 0 \\
			Schlaefli & 27 & 216 & \alphaopt{3} & 0 & \alphaopt{1} & 0 & \alphaopt{1} & 0 & \alphaopt{1} & 0 \\
			Coxeter & 28 & 42 & \alphaopt{12} & 0 & \alphaopt{7} & 0 & \alphaopt{4} & 0 & \alphaopt{1} & 0 \\
			DoubleStarSnark & 30 & 45 & \alphaopt{13} & 0 & \alphaopt{6} & 0 & \alphaopt{4} & 0 & \alphaopt{1} & 0 \\
			TutteCoxeter & 30 & 45 & \alphaopt{15} & 0 & \alphaopt{6} & 0 & \alphaopt{5} & 0 & \alphaopt{1} & 0 \\
			Dyck & 32 & 48 & \alphaopt{16} & 0 & \alphaopt{8} & 0 & \alphaopt{4} & 0 & \alphaopt{2} & 0 \\
			Wells & 32 & 80 & \alphaopt{10} & 0 & \alphaopt{2} & 0 & \alphaopt{2} & 0 & \alphaopt{1} & 0 \\
			Sylvester & 36 & 90 & \alphaopt{12} & 0 & \alphaopt{6} & 0 & \alphaopt{1} & 0 & \alphaopt{1} & 0 \\
			WienerAraya & 42 & 67 & \alphaopt{17} & 0 & \alphaopt{8} & 0 & \alphaopt{5} & 0 & \alphaopt{4} & 0 \\
			SzekeresSnark & 50 & 75 & \alphaopt{21} & 2 & \alphaopt{9} & 0 & \alphaopt{6} & 0 & \alphaopt{5} & 0 \\
			HoffmanSingleton & 50 & 175 & \alphaopt{15} & 0 & \alphaopt{1} & 0 & \alphaopt{1} & 0 & \alphaopt{1} & 0 \\
			Harborth & 52 & 104 & \alphaopt{18} & 0 & \alphaopt{10} & 0 & \alphaopt{6} & 0 & \alphaopt{4} & 0 \\
			EllinghamHorton54 & 54 & 81 & \alphaopt{27} & 0 & \alphaopt{11} & 0 & \alphaopt{8} & 0 & \alphaopt{4} & 0 \\
			Gray & 54 & 81 & \alphaopt{27} & 0 & \alphaopt{11} & 0 & \alphaopt{9} & 0 & \alphaopt{3} & 0 \\
			Klein3Regular & 56 & 84 & \alphaopt{23} & 3 & \alphaopt{12} & 0 & \alphaopt{7} & 0 & \alphaopt{4} & 0 \\
			SimsGewirtz & 56 & 280 & \alphaopt{16} & 0 & \alphaopt{1} & 0 & \alphaopt{1} & 0 & \alphaopt{1} & 0 \\
			Gosset & 56 & 756 & \alphaopt{4} & 0 & \alphaopt{2} & 0 & \alphaopt{1} & 0 & \alphaopt{1} & 0 \\
			Perkel & 57 & 171 & \alphaopt{19} & 0 & \alphaopt{5} & 0 & \alphaopt{1} & 0 & \alphaopt{1} & 0 \\
			BuckyBall & 60 & 90 & \alphaopt{24} & 7 & \alphaopt{12} & 0 & \alphaopt{7} & 0 & \alphaopt{6} & 0 \\
			Balaban10Cage & 70 & 105 & \alphaopt{35} & 0 & \alphaopt{17} & 0 & \alphaopt{9} & 1 & \alphaopt{5} & 0 \\
			HarriesWong & 70 & 105 & \alphaopt{35} & 0 & \alphaopt{17} & 0 & \alphaopt{9} & 1 & \alphaopt{5} & 0 \\
			Harries & 70 & 105 & \alphaopt{35} & 0 & \alphaopt{17} & 0 & \alphaopt{10} & 1 & \alphaopt{5} & 0 \\
			Meredith & 70 & 140 & \alphaopt{34} & 0 & \alphaopt{10} & 0 & \alphaopt{7} & 0 & \alphaopt{5} & 0 \\
			M22 & 77 & 616 & \alphaopt{21} & 0 & \alphaopt{1} & 0 & \alphaopt{1} & 0 & \alphaopt{1} & 0 \\
			EllinghamHorton78 & 78 & 117 & \alphaopt{39} & 0 & \alphaopt{18} & 1 & \alphaopt{11} & 0 & \alphaopt{6} & 0 \\
			BrouwerHaemers & 81 & 810 & \alphaopt{15} & 15 & \alphaopt{1} & 0 & \alphaopt{1} & 0 & \alphaopt{1} & 0 \\
			Foster & 90 & 135 & \alphaopt{45} & 1 & \alphaopt{21} & 2 & \alphaopt{15} & 1 & \alphaopt{5} & 5 \\
			Horton & 96 & 144 & \alphaopt{48} & 1 & \alphaopt{24} & 1 & \alphaopt{14} & 1 & \alphaopt{8} & 1 \\
			HallJanko & 100 & 1800 & \alphaopt{10} & 1 & \alphaopt{1} & 0 & \alphaopt{1} & 0 & \alphaopt{1} & 0 \\
			HigmanSims & 100 & 1100 & \alphaopt{22} & 1 & \alphaopt{1} & 0 & \alphaopt{1} & 0 & \alphaopt{1} & 0 \\
			Balaban11Cage & 112 & 168 & \alphaopt{52} & 15 & \alphaopt{24} & 21 & \alphaopt{16} & 6 & \alphaopt{9} & 4 \\
			Ljubljana & 112 & 168 & \alphaopt{56} & 2 & \alphaopt{26} & 15 & \alphaopt{17} & 5 & \alphaopt{8} & 13 \\
			Dejter & 112 & 336 & \alphaopt{56} & 2 & \alphaopt{16} & 2 & \alphaopt{8} & 3 & \alphaopt{2} & 5 \\
			Cell600 & 120 & 720 & \alphaopt{24} & 36 & \alphaopt{8} & 4 & \alphaopt{3} & 6 & \alphaopt{2} & 0 \\
			Tutte12Cage & 126 & 189 & \alphaopt{63} & 5 & \alphaopt{27} & 58 & \alphaopt{21} & 5 & \alphaopt{9} & 6 \\
			Cameron & 231 & 3465 & \alphaopt{21} & 22 & \alphaopt{1} & 0 & \alphaopt{1} & 0 & \alphaopt{1} & 0 \\
			Cell120 & 600 & 1200 & \alphanon{212} & 15280 & \alphaopt{120} & 645 & \alphanon{46} & 19091 & \alphanon{28} & 11832 \\
			
		\end{longtable}
	\end{center}
	
	\scriptsize
	\setlength{\tabcolsep}{3pt}
	\renewcommand{\arraystretch}{1.08}
	
	\begin{center}
		\begin{longtable}{l r r
				!{\;\vrule\;}
				r r !{\;\vrule\;} r r !{\;\vrule\;} r r !{\;\vrule\;} r r}
			\caption{Coding-theoretic instances, ordered by $|V|$. Runtimes are in seconds, rounded down. An asterisk indicates solver-reported optimality; $^\dagger$ identifies a result from a single-family configuration, with that run's time. Other entries use Both. Status and lower-bound conventions are given in Section~\ref{sec:benchmark}.}\label{tab:numbered-instances}\\
			\toprule
			Instance & $|V|$ & $|E|$ &
			$\alpha_1$ & Runtime &
			$\alpha_2$ & Runtime &
			$\alpha_3$ & Runtime &
			$\alpha_4$ & Runtime \\
			\midrule
			\endfirsthead
			
			\toprule
			Instance & $|V|$ & $|E|$ &
			$\alpha_1$ & Runtime &
			$\alpha_2$ & Runtime &
			$\alpha_3$ & Runtime &
			$\alpha_4$ & Runtime \\
			\midrule
			\endhead
			
			\bottomrule
			\endlastfoot
			
			\texttt{1tc.8}     & 8    & 6      & \alphaopt{4}   & 0    & \alphaopt{4}   & 0    & \alphaopt{4}   & 0     & \alphaopt{4}  & 0 \\
			
			\texttt{1tc.16}    & 16   & 22     & \alphaopt{8}   & 0    & \alphaopt{5}   & 0    & \alphaopt{5}   & 0     & \alphaopt{5}  & 0 \\
			
			\texttt{1tc.32}    & 32   & 68     & \alphaopt{12}  & 0    & \alphaopt{8}   & 0    & \alphaopt{6}   & 0     & \alphaopt{6}  & 0 \\
			
			\texttt{1dc.64}    & 64   & 543    & \alphaopt{10}  & 0    & \alphaopt{4}   & 0    & \alphaopt{2}   & 0     & \alphaopt{2}  & 0 \\
			\texttt{1et.64}    & 64   & 264    & \alphaopt{18}  & 0    & \alphaopt{8}   & 0    & \alphaopt{7}   & 0     & \alphaopt{7}  & 0 \\
			\texttt{1tc.64}    & 64   & 192    & \alphaopt{20}  & 0    & \alphaopt{12}  & 0    & \alphaopt{10}  & 0     & \alphaopt{8}  & 0 \\
			
			\texttt{1dc.128}   & 128  & 1471   & \alphaopt{16}  & 4    & \alphaopt{5}   & 7    & \alphaopt{2}   & 5     & \alphaopt{2}  & 1 \\
			\texttt{1et.128}   & 128  & 672    & \alphaopt{28}  & 2    & \alphaopt{14}  & 3    & \alphaopt{8}   & 2     & \alphaopt{8}  & 2 \\
			\texttt{1tc.128}   & 128  & 512    & \alphaopt{38}  & 2    & \alphaopt{20}  & 2    & \alphaopt{12}  & 2     & \alphaopt{12} & 3 \\
			\texttt{1zc.128}   & 128  & 1120   & \alphaopt{18}  & 20   & \alphaopt{4}   & 7    & \alphaopt{2}   & 6     & \alphaopt{2}  & 1 \\
			\texttt{2dc.128}   & 128  & 5173   & \alphaopt{5}   & 8    & \alphaopt{2}   & 1    & \alphaopt{2}   & 0     & \alphaopt{1}  & 0 \\
			
			\texttt{1dc.256}   & 256  & 3839   & \alphaopt{30}  & 88   & \alphaopt{7}   & 103  & \alphaopt{4}   & 144   & \alphaopt{2}  & 55 \\
			\texttt{1et.256}   & 256  & 1664   & \alphaopt{50}$^\dagger$  & 2059& \alphaopt{22}  & 39   & \alphaopt{12}  & 46    & \alphaopt{9}  & 46 \\
			\texttt{1tc.256}   & 256  & 1312   & \alphaopt{63}  & 23   & \alphaopt{32}  & 25   & \alphaopt{20}  & 34    & \alphaopt{14} & 36 \\
			\texttt{1zc.256}   & 256  & 2816   & \alphaopt{36}  & 51   & \alphaopt{7}   & 281  & \alphaopt{4}   & 150   & \alphaopt{2}  & 54 \\
			\texttt{2dc.256}   & 256  & 17183  & \alphaopt{7}   & 117  & \alphaopt{2}   & 55   & \alphaopt{2}   & 0     & \alphaopt{1}  & 0 \\
			
			\texttt{1dc.512}   & 512  & 9727   & \alphanon{52}  & 13631& \alphaopt{11}  & 7699 & \alphaopt{5}   & 2707  & \alphaopt{2}  & 1720 \\
			\texttt{1et.512}   & 512  & 4032   & \alphanon{100} & 11515& \alphaopt{30}  & 3381 & \alphaopt{22}  & 694   & \alphaopt{16} & 708 \\
			\texttt{1tc.512}   & 512  & 3264   & \alphaopt{110} & 3039 & \alphaopt{48}  & 433  & \alphaopt{30}  & 498   & \alphaopt{22} & 593 \\
			\texttt{1zc.512}   & 512  & 6912   & \alphanon{62}  & 11304& \alphanon{12}  & 11863& \alphaopt{4}   & 3466  & \alphaopt{2}  & 1728 \\
			\texttt{2dc.512}   & 512  & 54895  & \alphaopt{11}  & 7463 & \alphaopt{2}   & 1753 & \alphaopt{2}   & 59    & \alphaopt{2}  & 1 \\
			
			\texttt{1dc.1024}  & 1024 & 24063  & \alphanon{80}  & 23894& \alphanon{16}  & 14446& \alphaopt{6}   & 353   & \alphaopt{4}  & 5 \\
			\texttt{1et.1024}  & 1024 & 9600   & \alphanon{169} & 14086& \alphanon{58}  & 11299& \alphaopt{32}  & 356   & \alphaopt{21} & 88 \\
			\texttt{1tc.1024}  & 1024 & 7936   & \alphanon{195} & 14460& \alphanon{78}  & 11156& \alphaopt{51}  & 205   & \alphaopt{33} & 121 \\
			\texttt{1zc.1024}  & 1024 & 16640  & \alphanon{107} & 24863& \alphanon{18}  & 10982& \alphanon{6}   & 10820 & \alphaopt{4}  & 6 \\
			\texttt{2dc.1024}  & 1024 & 169162 & \alphanon{16}  & 14478& \alphaopt{4}   & 6    & \alphaopt{2}   & 4     & \alphaopt{2}  & 4 \\
			
			\texttt{1dc.2048}  & 2048 & 58367  & \alphanon{121} & 18585& \alphanon{22}  & 11307& \alphanon{8}   & 10913 & \alphaopt{5}  & 62 \\
			\texttt{1et.2048}  & 2048 & 22528  & \alphanon{291} & 29378& \alphanon{106} & 29393& \alphanon{44}  & 11915 & \alphaopt{32} & 1277 \\
			\texttt{1tc.2048}  & 2048 & 18944  & \alphanon{344} & 46082& \alphanon{130} & 20170& \alphanon{74}  & 13568 & \alphaopt{48} & 3779 \\
			\texttt{1zc.2048}  & 2048 & 39424  & \alphanon{173}  & 79132  & \alphanon{29}  & 22804  & \alphanon{8}   & 10963 & \alphanon{4}  & 10867 \\
			\texttt{2dc.2048}  & 2048 & 504451 & \alphanon{22}  & 11310& \alphaopt{5}   & 70   & \alphaopt{2}   & 44    & \alphaopt{2}  & 43 \\
			
			\texttt{1zc.4096}  & 4096 & 92160  & \alphanon{139} & 31875& \alphanon{12}  & 12277& \alphanon{12}  & 12523 & \alphanon{4}  & 11326 \\
			
		\end{longtable}
	\end{center}
	
	\begin{landscape}
		\begin{table}[ht]
			\centering
			\footnotesize
			\setlength{\tabcolsep}{3.5pt}
			\renewcommand{\arraystretch}{0.97}
			
			\caption{Sum-rank-metric instances, identified by $t,q,n,m$ as in Section~\ref{sec:benchmark}. The symbols $n,m$ denote tuples here. Runtimes are in seconds, rounded down. An asterisk indicates solver-reported optimality; $^\dagger$ identifies a single-family result and its runtime. The symbol $^{\S}$ identifies an unstarred entry whose value is known to be exact from~\cite{abiad2023eigenvalue}.}\label{tab:coding-results}
			
			\begin{tabular}{r r l l r r
					!{\;\vrule\;}
					r r !{\;\vrule\;} r r !{\;\vrule\;} r r !{\;\vrule\;} r r}
				\toprule
				$t$ & $q$ & $n$ & $m$ & $|V|$ & $|E|$ &
				$\alpha_1$ & Runtime &
				$\alpha_2$ & Runtime &
				$\alpha_3$ & Runtime &
				$\alpha_4$ & Runtime \\
				\midrule
				
				3  & 2 & (1,1,1) & (2,1,1) & 16   & 40     & \alphaopt{4}   & 0     & \alphaopt{2}   & 0     & \alphaopt{1}  & 0     & \alphaopt{1}  & 0 \\
				
				2  & 2 & (2,1) & (2,1) & 32   & 160    & \alphaopt{8}   & 0     & \alphaopt{2}   & 0     & \alphaopt{1}  & 0     & \alphaopt{1}  & 0 \\
				4  & 2 & (1,1,1,1) & (2,1,1,1) & 32   & 96     & \alphaopt{8}   & 0     & \alphaopt{4}   & 0     & \alphaopt{2}  & 0     & \alphaopt{1}  & 0 \\
				
				5  & 2 & (1,1,1,1,1) & (2,1,1,1,1) & 64   & 224    & \alphaopt{16}  & 0     & \alphaopt{8}   & 0     & \alphaopt{2}  & 0     & \alphaopt{2}  & 0 \\
				2  & 2 & (2,1) & (2,2) & 64   & 384    & \alphaopt{16}  & 0     & \alphaopt{4}   & 0     & \alphaopt{1}  & 0     & \alphaopt{1}  & 0 \\
				
				6  & 2 & (1,1,1,1,1,1) & (2,1,1,1,1,1) & 128  & 512    & \alphaopt{32}  & 2     & \alphaopt{9}   & 114   & \alphaopt{4}  & 7     & \alphaopt{2}  & 6 \\
				
				7  & 2 & (1,1,1,1,1,1,1) & (2,1,1,1,1,1,1) & 256  & 1152   & \alphaopt{64}  & 20    & \alphanon{18}  & 10839 & \alphaopt{8}  & 86    & \alphaopt{4}  & 145 \\
				4  & 2 & (2,1,1,1) & (2,2,1,1) & 256  & 1792   & \alphaopt{64}  & 51    & \alphaopt{16}  & 71    & \alphaopt{4}  & 134   & \alphaopt{2}  & 53 \\
				3  & 2 & (2,1,1) & (2,2,2) & 256  & 1920   & \alphaopt{64}  & 42    & \alphaopt{16}  & 80    & \alphaopt{4}  & 179   & \alphaopt{1}  & 0 \\
				2  & 2 & (2,2) & (2,2) & 256  & 2304   & \alphaopt{64}  & 52    & \alphaopt{9}   & 493   & \alphaopt{4}  & 176   & \alphaopt{1}  & 0 \\
				
				2  & 2 & (2,1) & (3,3) & 512  & 7168   & \alphaopt{64}  & 691   & \alphaopt{8}   & 1444  & \alphaopt{1}  & 0     & \alphaopt{1}  & 0 \\
				3  & 2 & (2,2,1) & (2,2,1) & 512  & 4864   & \alphaopt{128} & 1131  & \alphanon{18}$^{\S}$ & 11712 & \alphanon{5}  & 13002 & \alphaopt{2}  & 1736 \\
				4  & 2 & (2,1,1,1) & (2,2,2,1) & 512  & 4096   & \alphaopt{128} & 1143  & \alphanon{18}  & 11502 & \alphaopt{5}  & 2106  & \alphaopt{2}  & 1759 \\
				5  & 2 & (2,1,1,1,1) & (2,2,1,1,1) & 512  & 3840   & \alphaopt{128} & 684   & \alphaopt{32}  & 755   & \alphaopt{8}  & 1694  & \alphaopt{4}  & 2808 \\
				6  & 2 & (2,1,1,1,1,1) & (2,1,1,1,1,1) & 512  & 3584   & \alphaopt{128} & 556   & \alphaopt{32}  & 815   & \alphanon{9}$^{\S}$ & 12339 & \alphanon{4}  & 13528 \\
				8  & 2 & (1,1,1,1,1,1,1,1) & (2,1,1,1,1,1,1,1) & 512  & 2560   & \alphaopt{128} & 334   & \alphanon{32}  & 11292 & \alphaopt{16} & 1122  & \alphanon{5}  & 12941 \\
				
				2  & 3 & (2,1) & (2,2) & 729  & 14580  & \alphaopt{81}  & 3005  & \alphaopt{9}   & 15    & \alphaopt{1}  & 0     & \alphaopt{1}  & 0 \\
				
				7  & 2 & (2,1,1,1,1,1,1) & (2,1,1,1,1,1,1) & 1024 & 7680  & \alphaopt{256} & 8113  & \alphaopt{64}  & 10867 & \alphanon{18} & 11057 & \alphaopt{6}$^\dagger$  & 10817 \\
				
				8  & 2 & (2,1,1,1,1,1,1,1) & (2,1,1,1,1,1,1,1) & 2048 & 16384 & \alphaopt{512}$^\dagger$ & 150128& \alphanon{70}  & 24321 & \alphanon{31} & 25759 & \alphanon{9}  & 10955 \\
				
				9  & 2 & (2,1,1,1,1,1,1,1,1) & (2,1,1,1,1,1,1,1,1) & 4096 & 34816 & \alphanon{203} & 49495 & \alphanon{61}  & 19943 & \alphanon{12} & 12264 & \alphanon{12} & 20984 \\
				
				10 & 2 & (2,1,1,1,1,1,1,1,1,1) & (2,1,1,1,1,1,1,1,1,1) & 8192 & 73728  & \alphanon{159} & 137044& \alphanon{15}  & 23663 & \alphanon{13} & 22886 & \alphanon{10} & 20145 \\
				6  & 2 & (2,1,1,1,1,1) & (2,2,2,2,2,1) & 8192 & 90112  & \alphanon{362} & 295740& \alphanon{11}  & 22049 & \alphanon{11}  & 19873 & \alphanon{11} & 19120 \\
				6  & 2 & (2,2,1,1,1,1) & (2,2,2,1,1,1) & 8192 & 98304  & \alphanon{460} & 334121& \alphanon{16}  & 25514 & \alphanon{15} & 22212 & \alphanon{15} & 26674 \\
				5  & 2 & (2,2,1,1,1) & (2,2,2,2,1) & 8192 & 102400 & \alphanon{416} & 365426& \alphanon{10}  & 20648 & \alphanon{10} & 20303 & \alphanon{10} & 17944 \\
				4  & 2 & (2,2,2,1) & (2,2,2,1) & 8192 & 114688 & \alphanon{96}  & 88811 & \alphanon{10}  & 20876 & \alphanon{10} & 19995 & \alphanon{9}  & 18350 \\
				
				\bottomrule
			\end{tabular}
			
		\end{table}
	\end{landscape}

	\endgroup

	\subsection{Runtime profile across subproblems}\label{sec:effectQ}

Table~\ref{tab:perq-1et512} gives the recorded target-by-target profile for \texttt{1et.512}. It illustrates the distinction between finding a clique at a feasible target and excluding the next target. For an interrupted run, the largest clique found remains a lower bound rather than a determined value of $\alpha_k$.

\begin{table}[!tbp]
		\centering
		\scriptsize
		\setlength{\tabcolsep}{3pt}
		\renewcommand{\arraystretch}{1.0}
		\caption{Per-target runtimes for \texttt{1et.512} ($|V|=512$, $|E|=4032$). A range denotes targets with similar runtimes; only selected targets are shown for $k=1$. The $q=1$ time includes preprocessing. For $k=2,3,4$, the final target has solver-reported nonexistence of a $q$-clique, giving the starred value $q-1$ under Section~\ref{sec:benchmark}. For $k=1$, target $101$ times out (t.o.) and the reported lower bound is $100$. The counters \#cc and \#sc are cumulative counts of lifted rank inequalities~\eqref{cvi} and spectral rank inequalities~\eqref{sb-cut}, respectively; the broader counters in Appendix~\ref{app:ablation} also include the other tested cut families.}\label{tab:perq-1et512}
		
		\medskip
		
		\begin{tabular}{r r r r !{\;\vrule\;} r r r r}
			\toprule
			\multicolumn{4}{c!{\;\vrule\;}}{$k=1$ ($\alpha_1 \geq 100$)} & \multicolumn{4}{c}{$k=2$ ($\alpha_2 = 30^\ast$)} \\
			\cmidrule(lr){1-4}\cmidrule(lr){5-8}
			$q$ & Time (s) & \#cc & \#sc & $q$ & Time (s) & \#cc & \#sc \\
			\midrule
			1       & 102  & 0    & 0   & 1       & 408  & 0     & 0 \\
			2--53   & 0.6  & 0    & 0   & 2--26   & 0.6  & 0     & 0 \\
			54      & 0.7  & 50   & 1   & 27      & 0.9  & 138   & 5 \\
			55      & 1.0  & 100  & 2   & 28      & 1.4  & 256   & 10 \\
			60      & 1.5  & 372  & 6   & 29      & 7    & 1609 & 45 \\
			70      & 1.9  & 932  & 16  & 30      & 10   & 2555 & 72 \\
			80      & 2.2  & 1706 & 27 & 31  & 2937 & 3355 & 26744 \\
			90      & 4.1  & 4558 & 44 & & & & \\
			95      & 16   & 4966 & 46 & & & & \\
			99      & 202  & 5363 & 64 & & & & \\
			100     & 196  & 5436 & 65 & & & & \\
			101     & t.o. & 5479 & 66 & & & & \\
			\midrule
			\multicolumn{4}{c!{\;\vrule\;}}{$k=3$ ($\alpha_3 = 22^\ast$)} & \multicolumn{4}{c}{$k=4$ ($\alpha_4 = 16^\ast$)} \\
			\cmidrule(lr){1-4}\cmidrule(lr){5-8}
			$q$ & Time (s) & \#cc & \#sc & $q$ & Time (s) & \#cc & \#sc \\
			\midrule
			1       & 674  & 0    & 0   & 1       & 696  & 0     & 0 \\
			2--18   & 0.6  & 0    & 0   & 2--13   & 0.6  & 0     & 0 \\
			19      & 0.7  & 15   & 1   & 14      & 0.7  & 10    & 1 \\
			20      & 0.8  & 15   & 1   & 15      & 0.7  & 10    & 1 \\
			21      & 0.7  & 15   & 1   & 16      & 1.6  & 58    & 7 \\
			22      & 2.6  & 283  & 11  & 17      & 2.1  & 203   & 19 \\
			23      & 4.2  & 576  & 26  & & & & \\
			\bottomrule
		\end{tabular}
	\end{table}

The first target, $q=1$, includes construction of the auxiliary graph and the shortest-path inequalities. Its recorded time ranges from $102$ to $696$ seconds. Subsequent targets reuse this preprocessing. In the displayed ranges $q=2$--$53$, $2$--$26$, $2$--$18$, and $2$--$13$ for $k=1,2,3,4$, respectively, each target takes about $0.6$ seconds. These timings alone do not identify whether a clique was obtained by an integral LP solution, a primal heuristic, or branching.

For $k=2$, the successful targets through $q=30$ take approximately $440$ seconds in total, while the reported nonexistence conclusion at $q=31$ takes $2937$ seconds, about $87\%$ of the total runtime. The cumulative number of spectral rank cuts increases from $72$ to $26744$ at that target, so $26672$ are added during the final search. At the end of that search, the cumulative count exceeds the $3355$ lifted rank inequalities shown in the table. These two counters exclude other combinatorial and non-selected-vertex cuts, which are included in Appendix~\ref{app:ablation}.

For $k=1$, target $101$ reaches the time limit without finding a clique of that size or excluding it, leaving the reported lower bound $100$. This behavior motivates upper bounds that can close the search at a successful target. The recorded spectral bounds for $k=2,3,4$ are $166$, $86$, and $32$, respectively, so they do not close the runs at the starred values $30$, $22$, and $16$. A valid upper bound of $30$ for the $k=2$ case would have removed the final target-$31$ search. This is a consequence of the stopping rule, not a measured result for an alternative bound.

\subsection{Effect of the tested configurations}\label{sec:effectCuts}

The ablation uses three configurations. Both enables the combinatorial and spectral cuts and the spectral upper bound $q_{\max}$. Comb. enables only the combinatorial cuts and uses $q_{\max}=|V|$. Spec. enables the spectral cuts and the spectral upper bound, with no combinatorial cuts. All three use the same sequential framework and graph reduction. The counts in Table~\ref{tab:configuration-totals} are obtained from Appendix~\ref{app:ablation}; a success means solver-reported optimality under Section~\ref{sec:benchmark}.

\begin{table}[tbp]
\centering
\caption{Solver-reported successes of the tested configurations. The union counts combinations successful in at least one configuration and does not describe a run with a common total-time budget.}\label{tab:configuration-totals}
\begin{tabular}{lrrrrr}
\toprule
Class & Combinations & Both & Comb. & Spec. & Union \\
\midrule
Named graphs & 316 & 313 & 312 & 296 & 313 \\
Coding-theoretic graphs & 128 & 95 & 94 & 55 & 96 \\
Sum-rank-metric graphs & 100 & 62 & 63 & 42 & 64 \\
\midrule
Total & 544 & 470 & 469 & 393 & 473 \\
\bottomrule
\end{tabular}
\end{table}

The total number of successes is similar for Both and Comb., while Spec. resolves substantially fewer combinations. Both has three successes that Comb. lacks, and Comb. has two that Both lacks. Spec. supplies one further success absent from Both. These observations support a complementary role for the configurations, while the overall difference between Both and Comb. is modest.

The instance-level times show larger contrasts. For BuckyBall at $k=1$, Both returns a starred value in $7$ seconds whereas Spec. times out. Klein3Regular takes $3$ seconds with Both and $5496$ seconds with Spec. For Cell120 at $k=2$, Both returns a starred value of $120$ in $645$ seconds, whereas Spec. returns the lower bound $81$. On \texttt{1dc.256} at $k=1$, Both returns a starred value of $30$ in $88$ seconds, while Spec. reaches only the lower bound $26$ before timeout.

Spectral information can help through either separation or termination at $q_{\max}$. On Cameron at $k=1$, the recorded bound $q_{\max}=21$ matches the returned clique size, allowing Both to finish in $22$ seconds; Comb. times out while searching beyond its clique of size $21$. For \texttt{1et.512} at $k=2$, the bound $166$ is far above the returned value $30$, and Both completes target $31$ while Comb. times out. Conversely, on \texttt{1et.256} at $k=1$, Comb. returns the starred value $50$ in $2060$ seconds under the appendix's rounding convention, whereas Both times out.

Because the tested configurations change both spectral separation and the global bound, they do not isolate the cost or benefit of spectral cuts alone. A controlled comparison would retain the same $q_{\max}$ and preprocessing while switching only spectral separation, and would record separation time by family. The $z$-scaled inequalities and corrected non-selected-vertex cuts also require separate experiments before computational gains can be attributed to those theoretical strengthenings.

\subsection{Comparison with the discretized formulation of Martins}\label{sec:martins}

We compare the recorded results with the discretized formulation (DF) of Martins~\cite{martins2010extended}, specifically model (14)--(15) together with constraints (9), (10), (12), and (13) of that paper. DF solves the clique problem on $\bar G$ directly, without an outer loop over target sizes. Its implementation uses the same C++ and CPLEX environment, hardware, single-thread setting, and zero requested optimality gap. It receives the degree bound $q\leq\Delta(\bar G)+1$ but none of the additional cuts or spectral bounds. Models with more than $2\cdot10^8$ nonzero coefficients were excluded by the memory threshold and are marked mem.

In addition to the three main classes, this comparison includes $80$ DIMACS graphs~\cite{johnson1996cliques}. Each supplied graph is treated as the input $G$ to \textsc{$k$-MIS}. Consequently, at $k=1$ the objective is its independence number and the clique solver acts on its complement. These are not timings for the customary maximum-clique objective on the supplied DIMACS graph.

\paragraph{Budgets and aggregation}
DF has a nominal three-hour limit per instance and $k$, whereas the sequential method has a three-hour limit per target. The recorded comparison includes selected DF runs with extended limits. Across all four classes, the displayed runtime of our method exceeds three hours on $113$ combinations; only three of these are starred, all on sum-rank instances whose DF model was excluded. Appendix~\ref{app:martins} contains $38$ entries combining a displayed runtime above three hours for our method with a DF timeout. The dagger-marked entry for \texttt{1et.256} at $k=1$ displays the successful Comb. runtime instead of the longer unresolved Both runtime. The timing summaries therefore compare the recorded runs and configurations, rather than performance under a single common total-time budget.

The ours column of Table~\ref{tab:martins} includes the three dagger-marked single-family successes. This aggregation is the union of observed success sets, with the displayed runtime taken from the identified successful run; it excludes the cost of trying the other configurations. For the three main classes, Both alone has $55$ successes absent from DF, while the union has $58$. The corresponding DF-only count is one. These distinctions are necessary when interpreting the summary as evidence about an individual algorithm or a portfolio.

Table~\ref{tab:martins} retains the reported shifted geometric mean runtime on each common success set. For runtimes $t_1,\ldots,t_s$ measured in seconds, the statistic with a one-second shift is
\[
 \exp\!\left(\frac{1}{s}\sum_{j=1}^{s}\log(1+t_j)\right)-1.
\]
The times for our method were rounded down to whole seconds in these summaries. Small means and ratios involving them are sensitive to that rounding; the available aggregates do not support precise subsecond speed comparisons. The per-instance values and status counts are given in Appendix~\ref{app:martins}.

\paragraph{Observed outcomes}
On the named graphs, the reported results include $313$ successes for our configurations and $304$ for DF. DF has ten timeouts and two excluded models, and no DF-only success. On the coding-theoretic and sum-rank classes, our configurations have $24$ and $25$ successes absent from DF, respectively. DF has $19$ and $27$ timeouts in these classes, together with $37$ and $33$ excluded models. Its sole main-class success absent from our configurations is the $512$-vertex sum-rank graph with $t=6$, $q=2$, $n=m=(2,1^5)$ at $k=4$, for which DF returns $4$ after $6581$ seconds.

On DIMACS at $k=1$, our method has $41$ successes and DF has $36$, with $33$ in common. Eight are exclusive to our method and three---\texttt{brock200\_1}, \texttt{brock200\_4}, and \texttt{sanr200\_0.7}---are exclusive to DF. For $k\geq2$, our method has $80$ successes at each value of $k$, and DF has $77$; the DF timeouts are \texttt{c-fat500-1}, \texttt{c-fat500-2}, and \texttt{c-fat500-5}. Whenever both methods report optimality, their objective values agree. This consistency is useful evidence about the computations, subject to the numerical-status convention already stated.

\begin{table}[!htbp]
		\centering
		\footnotesize
		\setlength{\tabcolsep}{4pt}
		\renewcommand{\arraystretch}{1.0}
		\caption{Comparison with the discretized formulation (DF) of Martins~\cite{martins2010extended}. For each family and $k$: number of combinations ($n$); solver-reported successes for the reported results of our configurations (ours) and DF; DF timeouts (t.o.) and models excluded by the nonzero threshold (mem); and overlap of the success sets. The ours column includes the three dagger-marked single-family results; it is not the fixed Both configuration. Mean runtimes use the common success set and a shift of $1\,\mathrm{s}$; the timing summaries retain whole-second rounding for our method. Section~\ref{sec:martins} explains the budget and aggregation limitations, and Section~\ref{sec:benchmark} defines solver-reported optimality.}\label{tab:martins}
		\begin{tabular}{l r r r r r r r r r r}
			\toprule
			Family / $k$ & $n$ & ours & DF & t.o. & mem & only ours & only DF & both & $\overline{t}_{\text{ours}}$ & $\overline{t}_{\text{DF}}$ \\
			\midrule
			\multicolumn{11}{l}{\textit{Named graphs}} \\
			\quad $k=1$ & 79 & 78 & 76 & 2 & 1 & 2 & 0 & 76 & 0.2 & 2.7 \\
			\quad $k=2$ & 79 & 79 & 75 & 3 & 1 & 4 & 0 & 75 & 0.1 & 1.7 \\
			\quad $k=3$ & 79 & 78 & 75 & 4 & 0 & 3 & 0 & 75 & 0.1 & 0.9 \\
			\quad $k=4$ & 79 & 78 & 78 & 1 & 0 & 0 & 0 & 78 & 0.1 & 0.7 \\
			\quad all & 316 & 313 & 304 & 10 & 2 & 9 & 0 & 304 & 0.1 & 1.4 \\
			\addlinespace
			\multicolumn{11}{l}{\textit{Coding-theoretic}} \\
			\quad $k=1$ & 32 & 18 & 13 & 8 & 11 & 5 & 0 & 13 & 2.7 & 62.1 \\
			\quad $k=2$ & 32 & 22 & 17 & 5 & 10 & 5 & 0 & 17 & 7.3 & 23.8 \\
			\quad $k=3$ & 32 & 26 & 20 & 3 & 9 & 6 & 0 & 20 & 6.8 & 11.9 \\
			\quad $k=4$ & 32 & 30 & 22 & 3 & 7 & 8 & 0 & 22 & 5.8 & 7.8 \\
			\quad all & 128 & 96 & 72 & 19 & 37 & 24 & 0 & 72 & 5.6 & 16.8 \\
			\addlinespace
			\multicolumn{11}{l}{\textit{Sum-rank-metric graphs}} \\
			\quad $k=1$ & 25 & 19 & 5 & 11 & 9 & 14 & 0 & 5 & 0.0 & 2.6 \\
			\quad $k=2$ & 25 & 14 & 7 & 10 & 8 & 7 & 0 & 7 & 3.8 & 9.9 \\
			\quad $k=3$ & 25 & 15 & 12 & 5 & 8 & 3 & 0 & 12 & 5.2 & 6.2 \\
			\quad $k=4$ & 25 & 16 & 16 & 1 & 8 & 1 & 1 & 15 & 8.5 & 3.7 \\
			\quad all & 100 & 64 & 40 & 27 & 33 & 25 & 1 & 39 & 4.5 & 5.0 \\
			\addlinespace
			\multicolumn{11}{l}{\textit{DIMACS}} \\
			\quad $k=1$ & 80 & 41 & 36 & 33 & 11 & 8 & 3 & 33 & 51.2 & 26.6 \\
			\quad $k=2$ & 80 & 80 & 77 & 3 & 0 & 3 & 0 & 77 & 0.3 & 0.4 \\
			\quad $k=3$ & 80 & 80 & 77 & 3 & 0 & 3 & 0 & 77 & 0.2 & 0.3 \\
			\quad $k=4$ & 80 & 80 & 77 & 3 & 0 & 3 & 0 & 77 & 0.1 & 0.2 \\
			\quad all & 320 & 281 & 267 & 42 & 11 & 17 & 3 & 264 & 0.9 & 0.9 \\
			\bottomrule
		\end{tabular}
	\end{table}

\section{Conclusion}\label{sec:conclusion}

We have developed a sequential branch-and-cut framework for maximum $k$-independent sets in which each target clique size gives a compact model with one binary variable per remaining vertex and a continuous activation variable. We have proved the exactness of this model, derived activation-scaled rank inequalities and their liftings, and strengthened the non-selected-vertex inequality. These cuts dominate their unscaled counterparts in the fixed-size relaxation. We use classical spectral bounds both for global stopping and for local rank certificates on complements of subgraphs of the current working graph. Evaluating the bounds on that graph is essential after edge peeling. Under valid bounds, valid cuts, and completed solves, Theorem~\ref{thm:framework-correctness} establishes exactness independently of the completeness of the separation heuristics.

The reported implementation gives $470$ solver-reported successes among $544$ main-benchmark combinations, and the union with the single-family runs gives $473$. The ablation identifies instances on which spectral information helps and others on which its additional work is unhelpful. It also shows that spectral separation must be distinguished from the stopping bound. The comparison with Martins' discretized formulation yields more reported successes for the sequential configurations on these benchmark classes, while also identifying DF-only successes and limitations caused by differing budgets and configuration aggregation.

Several developments follow directly from the analysis. The continuous-$z$ formulation without the redundant rows, $z$-scaled cuts, and corrected non-selected-vertex inequalities should be evaluated against the computational variant studied here. A proof-producing implementation requires certified spectral calculations or independent validation of every affected bound and cut, together with stored clique witnesses. Comparisons using a common total-time budget, identical global bounds where appropriate, unrounded runtimes, and specialised maximum-clique solvers would clarify the computational contribution of each ingredient. Polynomial spectral bounds obtained directly from $G$~\cite{ACFNS2022,abiad2024unified,abiad2026inertia}, improved primal heuristics, and implicit representations of the auxiliary graph provide further directions for reducing preprocessing and unsuccessful-target costs.

\subsection*{Acknowledgements}
Aida Abiad is supported by NWO (Dutch Research Council) through the grants VI.Vidi.213.085 and OCENW.KLEIN.475.

\paragraph{Use of artificial intelligence}
ChatGPT (OpenAI) was used to assist with mathematical proofreading, exposition, and LaTeX preparation.
The authors assume responsibility for all content.

\affl{Aida Abiad}{a.abiad.monge@tue.nl}{Department of Mathematics and Computer Science, Eindhoven University of Technology, The Netherlands \\
Department of Mathematics and Data Science, Vrije Universiteit Brussel, Belgium}

\affl{Julian Golak}{julian.golak@uni-hamburg.de}{University of Hamburg Business School, Institute of Operations Management, Germany}

\affl{Arne Schulz}{arne.schulz@uni-hamburg.de}{University of Hamburg Business School, Institute of Operations Management, Germany}

\clearpage
\appendix
\renewcommand{\thetable}{A\arabic{table}}\setcounter{table}{0}
\section{Per-instance comparison with the discretized formulation}\label{app:martins}

Table~\ref{tab:martins-full} gives the reported values and runtimes of our configurations and the DF formulation under Section~\ref{sec:martins}. The starred and dagger-marked entries use the conventions of Section~\ref{sec:benchmark}; the ours column incorporates the three identified single-family substitutions and is not a fixed-Both comparison. Unstarred entries are reported lower bounds. The table includes the DIMACS results, which are not repeated in the main instance tables. For the sum-rank instances, $x^j$ abbreviates $j$ consecutive entries equal to $x$ in a parameter tuple.

\begingroup
	\begin{landscape}
	\scriptsize
	\setlength{\tabcolsep}{2pt}
	\renewcommand{\arraystretch}{1.05}
	\begin{center}
\begin{longtable}{l r !{\,\vrule\,} r r r r !{\,\vrule\,} r r r r !{\,\vrule\,} r r r r !{\,\vrule\,} r r r r}
		\caption{Per-instance comparison with DF~\cite{martins2010extended}. For each $k$, columns give the reported value and runtime in seconds, rounded down. An asterisk denotes solver-reported optimality as defined in Section~\ref{sec:benchmark}. For our configurations, $^\dagger$ marks a single-family result and its time in place of the unresolved Both run. For DF, t.o. denotes a timeout under the recorded budget and the adjacent value is its incumbent; -- denotes a model excluded by the nonzero threshold. Budgets and aggregation are described in Section~\ref{sec:martins}.}\label{tab:martins-full}\\
		\toprule
		Instance & $|V|$ & \multicolumn{4}{c}{$k=1$} & \multicolumn{4}{c}{$k=2$} & \multicolumn{4}{c}{$k=3$} & \multicolumn{4}{c}{$k=4$} \\
		 & & $\alpha$ ours & $t$ & $\alpha$ DF & $t$ & $\alpha$ ours & $t$ & $\alpha$ DF & $t$ & $\alpha$ ours & $t$ & $\alpha$ DF & $t$ & $\alpha$ ours & $t$ & $\alpha$ DF & $t$ \\
		\midrule
		\endfirsthead
		\toprule
		Instance & $|V|$ & \multicolumn{4}{c}{$k=1$} & \multicolumn{4}{c}{$k=2$} & \multicolumn{4}{c}{$k=3$} & \multicolumn{4}{c}{$k=4$} \\
		 & & $\alpha$ ours & $t$ & $\alpha$ DF & $t$ & $\alpha$ ours & $t$ & $\alpha$ DF & $t$ & $\alpha$ ours & $t$ & $\alpha$ DF & $t$ & $\alpha$ ours & $t$ & $\alpha$ DF & $t$ \\
		\midrule
		\endhead
		\bottomrule
		\endlastfoot
		\multicolumn{18}{l}{\textit{Named graphs}} \\
		Tetrahedral & 4 & \alphaopt{1} & 0 & \alphaopt{1} & 0 & \alphaopt{1} & 0 & \alphaopt{1} & 0 & \alphaopt{1} & 0 & \alphaopt{1} & 0 & \alphaopt{1} & 0 & \alphaopt{1} & 0 \\
		Octahedral & 6 & \alphaopt{2} & 0 & \alphaopt{2} & 0 & \alphaopt{1} & 0 & \alphaopt{1} & 0 & \alphaopt{1} & 0 & \alphaopt{1} & 0 & \alphaopt{1} & 0 & \alphaopt{1} & 0 \\
		Thomsen & 6 & \alphaopt{3} & 0 & \alphaopt{3} & 0 & \alphaopt{1} & 0 & \alphaopt{1} & 0 & \alphaopt{1} & 0 & \alphaopt{1} & 0 & \alphaopt{1} & 0 & \alphaopt{1} & 0 \\
		MoserSpindle & 7 & \alphaopt{2} & 0 & \alphaopt{2} & 0 & \alphaopt{1} & 0 & \alphaopt{1} & 0 & \alphaopt{1} & 0 & \alphaopt{1} & 0 & \alphaopt{1} & 0 & \alphaopt{1} & 0 \\
		Hexahedral & 8 & \alphaopt{4} & 0 & \alphaopt{4} & 0 & \alphaopt{2} & 0 & \alphaopt{2} & 0 & \alphaopt{1} & 0 & \alphaopt{1} & 0 & \alphaopt{1} & 0 & \alphaopt{1} & 0 \\
		Wagner & 8 & \alphaopt{3} & 0 & \alphaopt{3} & 0 & \alphaopt{1} & 0 & \alphaopt{1} & 0 & \alphaopt{1} & 0 & \alphaopt{1} & 0 & \alphaopt{1} & 0 & \alphaopt{1} & 0 \\
		KrackhardtKite & 10 & \alphaopt{4} & 0 & \alphaopt{4} & 0 & \alphaopt{2} & 0 & \alphaopt{2} & 0 & \alphaopt{2} & 0 & \alphaopt{2} & 0 & \alphaopt{1} & 0 & \alphaopt{1} & 0 \\
\end{longtable}
	\end{center}
	\end{landscape}
	\endgroup

	\section{Ablation tables}\label{app:ablation}
\renewcommand{\thetable}{B\arabic{table}}\setcounter{table}{0}

Tables~\ref{tab:named-both-vs-comb}--\ref{tab:graph-both-vs-spec} report the individual runs underlying Section~\ref{sec:effectCuts}. Each group gives the returned value, runtime, and accumulated cut counts for Both and for Comb. or Spec. An asterisk denotes solver-reported optimality under Section~\ref{sec:benchmark}. Each value belongs to the indicated configuration and $k$; bounds are not transferred across $k$ or across configurations. Consequently, a single-family column can contain a larger lower bound than the corresponding main-table entry.

The counter cc includes anticlique, lifted rank, and shortest-path inequalities; sc includes spectral rank and the tested weaker non-selected-vertex inequalities~\eqref{nsv-legacy}. These are the unscaled inequalities used in the computational variant of Section~\ref{sec:solutionProceedure}. Counters accumulate over all targets. Runtimes are rounded to the nearest second and may differ by one second from the rounded-down values in the main tables and Appendix~\ref{app:martins}.

\begingroup
\fontsize{7}{8.2}\selectfont
\setlength{\tabcolsep}{2pt}
\renewcommand{\arraystretch}{1}
\setlength{\LTleft}{\fill}
\setlength{\LTright}{\fill}
\begin{landscape}
\clearpage
\begin{longtable}{@{}l !{\,\vrule\,} r r r r !{\,\vrule\,} r r r !{\,\vrule\,} r r r r !{\,\vrule\,} r r r @{}}
\caption{Named graph instances: both cut families (Both) versus combinatorial cuts only (Comb.). Panel (a): $k=1,2$.}\label{tab:named-both-vs-comb}\\
\toprule
\multicolumn{1}{c}{} & \multicolumn{7}{c}{$k=1$} & \multicolumn{7}{c}{$k=2$} \\
\cmidrule(lr){2-8}\cmidrule(lr){9-15}
\multicolumn{1}{c}{} & \multicolumn{4}{c}{Both} & \multicolumn{3}{c}{Comb.} & \multicolumn{4}{c}{Both} & \multicolumn{3}{c}{Comb.} \\
\cmidrule(lr){2-5}\cmidrule(lr){6-8}\cmidrule(lr){9-12}\cmidrule(lr){13-15}
Instance & $\alpha$ & time & cc & sc & $\alpha$ & time & cc & $\alpha$ & time & cc & sc & $\alpha$ & time & cc \\
\midrule
\endfirsthead
\caption[]{Named graph instances: both cut families (Both) versus combinatorial cuts only (Comb.). Panel (a), continued.}\\
\toprule
\multicolumn{1}{c}{} & \multicolumn{7}{c}{$k=1$} & \multicolumn{7}{c}{$k=2$} \\
\cmidrule(lr){2-8}\cmidrule(lr){9-15}
\multicolumn{1}{c}{} & \multicolumn{4}{c}{Both} & \multicolumn{3}{c}{Comb.} & \multicolumn{4}{c}{Both} & \multicolumn{3}{c}{Comb.} \\
\cmidrule(lr){2-5}\cmidrule(lr){6-8}\cmidrule(lr){9-12}\cmidrule(lr){13-15}
Instance & $\alpha$ & time & cc & sc & $\alpha$ & time & cc & $\alpha$ & time & cc & sc & $\alpha$ & time & cc \\
\midrule
\endhead
\midrule
\multicolumn{15}{r}{Continued on next page} \\
\endfoot
\bottomrule
\endlastfoot
Tetrahedral & \alphaopt{1} & 0 & 0 & 0 & \alphaopt{1} & 0 & 0 & \alphaopt{1} & 0 & 0 & 0 & \alphaopt{1} & 0 & 0 \\
Thomsen & \alphaopt{3} & 0 & 0 & 0 & \alphaopt{3} & 0 & 0 & \alphaopt{1} & 0 & 0 & 0 & \alphaopt{1} & 0 & 0 \\
Octahedral & \alphaopt{2} & 0 & 0 & 0 & \alphaopt{2} & 0 & 0 & \alphaopt{1} & 0 & 0 & 0 & \alphaopt{1} & 0 & 0 \\
MoserSpindle & \alphaopt{2} & 0 & 11 & 0 & \alphaopt{2} & 0 & 11 & \alphaopt{1} & 0 & 0 & 0 & \alphaopt{1} & 0 & 0 \\
Hexahedral & \alphaopt{4} & 0 & 12 & 0 & \alphaopt{4} & 0 & 12 & \alphaopt{2} & 0 & 0 & 0 & \alphaopt{2} & 0 & 0 \\
Wagner & \alphaopt{3} & 0 & 12 & 0 & \alphaopt{3} & 0 & 12 & \alphaopt{1} & 0 & 0 & 0 & \alphaopt{1} & 0 & 0 \\
Petersen & \alphaopt{4} & 0 & 15 & 0 & \alphaopt{4} & 0 & 23 & \alphaopt{1} & 0 & 0 & 0 & \alphaopt{1} & 0 & 0 \\
KrackhardtKite & \alphaopt{4} & 0 & 20 & 1 & \alphaopt{4} & 0 & 20 & \alphaopt{2} & 0 & 22 & 0 & \alphaopt{2} & 0 & 22 \\
\end{longtable}
\clearpage
\addtocounter{table}{-1}
\begin{longtable}{@{}l !{\,\vrule\,} r r r r !{\,\vrule\,} r r r !{\,\vrule\,} r r r r !{\,\vrule\,} r r r @{}}
\caption{Named graph instances: both cut families (Both) versus combinatorial cuts only (Comb.). Panel (b): $k=3,4$. Continued.}\\
\toprule
\multicolumn{1}{c}{} & \multicolumn{7}{c}{$k=3$} & \multicolumn{7}{c}{$k=4$} \\
\cmidrule(lr){2-8}\cmidrule(lr){9-15}
\multicolumn{1}{c}{} & \multicolumn{4}{c}{Both} & \multicolumn{3}{c}{Comb.} & \multicolumn{4}{c}{Both} & \multicolumn{3}{c}{Comb.} \\
\cmidrule(lr){2-5}\cmidrule(lr){6-8}\cmidrule(lr){9-12}\cmidrule(lr){13-15}
Instance & $\alpha$ & time & cc & sc & $\alpha$ & time & cc & $\alpha$ & time & cc & sc & $\alpha$ & time & cc \\
\midrule
\endfirsthead
\caption[]{Named graph instances: both cut families (Both) versus combinatorial cuts only (Comb.). Panel (b), continued.}\\
\toprule
\multicolumn{1}{c}{} & \multicolumn{7}{c}{$k=3$} & \multicolumn{7}{c}{$k=4$} \\
\cmidrule(lr){2-8}\cmidrule(lr){9-15}
\multicolumn{1}{c}{} & \multicolumn{4}{c}{Both} & \multicolumn{3}{c}{Comb.} & \multicolumn{4}{c}{Both} & \multicolumn{3}{c}{Comb.} \\
\cmidrule(lr){2-5}\cmidrule(lr){6-8}\cmidrule(lr){9-12}\cmidrule(lr){13-15}
Instance & $\alpha$ & time & cc & sc & $\alpha$ & time & cc & $\alpha$ & time & cc & sc & $\alpha$ & time & cc \\
\midrule
\endhead
\midrule
\multicolumn{15}{r}{Continued on next page} \\
\endfoot
\bottomrule
\endlastfoot
Tetrahedral & \alphaopt{1} & 0 & 0 & 0 & \alphaopt{1} & 0 & 0 & \alphaopt{1} & 0 & 0 & 0 & \alphaopt{1} & 0 & 0 \\
Thomsen & \alphaopt{1} & 0 & 0 & 0 & \alphaopt{1} & 0 & 0 & \alphaopt{1} & 0 & 0 & 0 & \alphaopt{1} & 0 & 0 \\
Octahedral & \alphaopt{1} & 0 & 0 & 0 & \alphaopt{1} & 0 & 0 & \alphaopt{1} & 0 & 0 & 0 & \alphaopt{1} & 0 & 0 \\
MoserSpindle & \alphaopt{1} & 0 & 0 & 0 & \alphaopt{1} & 0 & 0 & \alphaopt{1} & 0 & 0 & 0 & \alphaopt{1} & 0 & 0 \\
Hexahedral & \alphaopt{1} & 0 & 0 & 0 & \alphaopt{1} & 0 & 0 & \alphaopt{1} & 0 & 0 & 0 & \alphaopt{1} & 0 & 0 \\
Wagner & \alphaopt{1} & 0 & 0 & 0 & \alphaopt{1} & 0 & 0 & \alphaopt{1} & 0 & 0 & 0 & \alphaopt{1} & 0 & 0 \\
Petersen & \alphaopt{1} & 0 & 0 & 0 & \alphaopt{1} & 0 & 0 & \alphaopt{1} & 0 & 0 & 0 & \alphaopt{1} & 0 & 0 \\
KrackhardtKite & \alphaopt{2} & 0 & 10 & 0 & \alphaopt{2} & 0 & 10 & \alphaopt{1} & 0 & 0 & 0 & \alphaopt{1} & 0 & 0 \\
\end{longtable}
\clearpage
\begin{longtable}{@{}l !{\,\vrule\,} r r r r !{\,\vrule\,} r r r !{\,\vrule\,} r r r r !{\,\vrule\,} r r r @{}}
\caption{Coding-theoretic instances: both cut families (Both) versus combinatorial cuts only (Comb.). Panel (a): $k=1,2$.}\label{tab:numeric-both-vs-comb}\\
\toprule
\multicolumn{1}{c}{} & \multicolumn{7}{c}{$k=1$} & \multicolumn{7}{c}{$k=2$} \\
\cmidrule(lr){2-8}\cmidrule(lr){9-15}
\multicolumn{1}{c}{} & \multicolumn{4}{c}{Both} & \multicolumn{3}{c}{Comb.} & \multicolumn{4}{c}{Both} & \multicolumn{3}{c}{Comb.} \\
\cmidrule(lr){2-5}\cmidrule(lr){6-8}\cmidrule(lr){9-12}\cmidrule(lr){13-15}
Instance & $\alpha$ & time & cc & sc & $\alpha$ & time & cc & $\alpha$ & time & cc & sc & $\alpha$ & time & cc \\
\midrule
\endfirsthead
\caption[]{Coding-theoretic instances: both cut families (Both) versus combinatorial cuts only (Comb.). Panel (a), continued.}\\
\toprule
\multicolumn{1}{c}{} & \multicolumn{7}{c}{$k=1$} & \multicolumn{7}{c}{$k=2$} \\
\cmidrule(lr){2-8}\cmidrule(lr){9-15}
\multicolumn{1}{c}{} & \multicolumn{4}{c}{Both} & \multicolumn{3}{c}{Comb.} & \multicolumn{4}{c}{Both} & \multicolumn{3}{c}{Comb.} \\
\cmidrule(lr){2-5}\cmidrule(lr){6-8}\cmidrule(lr){9-12}\cmidrule(lr){13-15}
Instance & $\alpha$ & time & cc & sc & $\alpha$ & time & cc & $\alpha$ & time & cc & sc & $\alpha$ & time & cc \\
\midrule
\endhead
\midrule
\multicolumn{15}{r}{Continued on next page} \\
\endfoot
\bottomrule
\endlastfoot
1tc.8 & \alphaopt{4} & 0 & 6 & 0 & \alphaopt{4} & 0 & 6 & \alphaopt{4} & 0 & 6 & 0 & \alphaopt{4} & 0 & 6 \\
1tc.16 & \alphaopt{8} & 0 & 35 & 5 & \alphaopt{8} & 0 & 32 & \alphaopt{5} & 0 & 27 & 0 & \alphaopt{5} & 0 & 31 \\
1tc.32 & \alphaopt{12} & 0 & 98 & 3 & \alphaopt{12} & 0 & 119 & \alphaopt{8} & 0 & 125 & 4 & \alphaopt{8} & 0 & 125 \\
1dc.64 & \alphaopt{10} & 0 & 659 & 30 & \alphaopt{10} & 0 & 680 & \alphaopt{4} & 0 & 1398 & 2 & \alphaopt{4} & 0 & 1398 \\
1et.64 & \alphaopt{18} & 0 & 635 & 30 & \alphaopt{18} & 0 & 674 & \alphaopt{8} & 0 & 452 & 6 & \alphaopt{8} & 0 & 452 \\
1tc.64 & \alphaopt{20} & 0 & 478 & 13 & \alphaopt{20} & 0 & 478 & \alphaopt{12} & 0 & 467 & 18 & \alphaopt{12} & 0 & 466 \\
1dc.128 & \alphaopt{16} & 4 & 2017 & 50 & \alphaopt{16} & 5 & 1986 & \alphaopt{5} & 7 & 5237 & 24 & \alphaopt{5} & 13 & 5230 \\
1et.128 & \alphaopt{28} & 3 & 1639 & 26 & \alphaopt{28} & 3 & 1520 & \alphaopt{14} & 3 & 1717 & 20 & \alphaopt{14} & 3 & 1723 \\
\end{longtable}
\clearpage
\addtocounter{table}{-1}
\begin{longtable}{@{}l !{\,\vrule\,} r r r r !{\,\vrule\,} r r r !{\,\vrule\,} r r r r !{\,\vrule\,} r r r @{}}
\caption{Coding-theoretic instances: both cut families (Both) versus combinatorial cuts only (Comb.). Panel (b): $k=3,4$. Continued.}\\
\toprule
\multicolumn{1}{c}{} & \multicolumn{7}{c}{$k=3$} & \multicolumn{7}{c}{$k=4$} \\
\cmidrule(lr){2-8}\cmidrule(lr){9-15}
\multicolumn{1}{c}{} & \multicolumn{4}{c}{Both} & \multicolumn{3}{c}{Comb.} & \multicolumn{4}{c}{Both} & \multicolumn{3}{c}{Comb.} \\
\cmidrule(lr){2-5}\cmidrule(lr){6-8}\cmidrule(lr){9-12}\cmidrule(lr){13-15}
Instance & $\alpha$ & time & cc & sc & $\alpha$ & time & cc & $\alpha$ & time & cc & sc & $\alpha$ & time & cc \\
\midrule
\endfirsthead
\caption[]{Coding-theoretic instances: both cut families (Both) versus combinatorial cuts only (Comb.). Panel (b), continued.}\\
\toprule
\multicolumn{1}{c}{} & \multicolumn{7}{c}{$k=3$} & \multicolumn{7}{c}{$k=4$} \\
\cmidrule(lr){2-8}\cmidrule(lr){9-15}
\multicolumn{1}{c}{} & \multicolumn{4}{c}{Both} & \multicolumn{3}{c}{Comb.} & \multicolumn{4}{c}{Both} & \multicolumn{3}{c}{Comb.} \\
\cmidrule(lr){2-5}\cmidrule(lr){6-8}\cmidrule(lr){9-12}\cmidrule(lr){13-15}
Instance & $\alpha$ & time & cc & sc & $\alpha$ & time & cc & $\alpha$ & time & cc & sc & $\alpha$ & time & cc \\
\midrule
\endhead
\midrule
\multicolumn{15}{r}{Continued on next page} \\
\endfoot
\bottomrule
\endlastfoot
1tc.8 & \alphaopt{4} & 0 & 6 & 0 & \alphaopt{4} & 0 & 6 & \alphaopt{4} & 0 & 6 & 0 & \alphaopt{4} & 0 & 6 \\
1tc.16 & \alphaopt{5} & 0 & 27 & 0 & \alphaopt{5} & 0 & 31 & \alphaopt{5} & 0 & 27 & 0 & \alphaopt{5} & 0 & 31 \\
1tc.32 & \alphaopt{6} & 0 & 110 & 0 & \alphaopt{6} & 0 & 120 & \alphaopt{6} & 0 & 110 & 0 & \alphaopt{6} & 0 & 120 \\
1dc.64 & \alphaopt{2} & 0 & 462 & 0 & \alphaopt{2} & 0 & 462 & \alphaopt{2} & 0 & 42 & 0 & \alphaopt{2} & 0 & 42 \\
1et.64 & \alphaopt{7} & 0 & 430 & 0 & \alphaopt{7} & 0 & 447 & \alphaopt{7} & 0 & 430 & 0 & \alphaopt{7} & 0 & 447 \\
1tc.64 & \alphaopt{10} & 0 & 460 & 6 & \alphaopt{10} & 0 & 460 & \alphaopt{8} & 0 & 435 & 1 & \alphaopt{8} & 0 & 448 \\
1dc.128 & \alphaopt{2} & 6 & 4032 & 0 & \alphaopt{2} & 9 & 4032 & \alphaopt{2} & 1 & 812 & 0 & \alphaopt{2} & 2 & 812 \\
1et.128 & \alphaopt{8} & 3 & 1652 & 0 & \alphaopt{8} & 4 & 1673 & \alphaopt{8} & 3 & 1652 & 0 & \alphaopt{8} & 4 & 1673 \\
\end{longtable}
\clearpage
\begin{longtable}{@{}r r l r !{\,\vrule\,} r r r r !{\,\vrule\,} r r r !{\,\vrule\,} r r r r !{\,\vrule\,} r r r @{}}
\caption{Sum-rank-metric graph instances: both cut families (Both) versus combinatorial cuts only (Comb.). Panel (a): $k=1,2$.}\label{tab:graph-both-vs-comb}\\
\toprule
\multicolumn{4}{c}{} & \multicolumn{7}{c}{$k=1$} & \multicolumn{7}{c}{$k=2$} \\
\cmidrule(lr){5-11}\cmidrule(lr){12-18}
\multicolumn{4}{c}{} & \multicolumn{4}{c}{Both} & \multicolumn{3}{c}{Comb.} & \multicolumn{4}{c}{Both} & \multicolumn{3}{c}{Comb.} \\
\cmidrule(lr){5-8}\cmidrule(lr){9-11}\cmidrule(lr){12-15}\cmidrule(lr){16-18}
$t$ & $q$ & $n$, $m$ & $|V|$ & $\alpha$ & time & cc & sc & $\alpha$ & time & cc & $\alpha$ & time & cc & sc & $\alpha$ & time & cc \\
\midrule
\endfirsthead
\caption[]{Sum-rank-metric graph instances: both cut families (Both) versus combinatorial cuts only (Comb.). Panel (a), continued.}\\
\toprule
\multicolumn{4}{c}{} & \multicolumn{7}{c}{$k=1$} & \multicolumn{7}{c}{$k=2$} \\
\cmidrule(lr){5-11}\cmidrule(lr){12-18}
\multicolumn{4}{c}{} & \multicolumn{4}{c}{Both} & \multicolumn{3}{c}{Comb.} & \multicolumn{4}{c}{Both} & \multicolumn{3}{c}{Comb.} \\
\cmidrule(lr){5-8}\cmidrule(lr){9-11}\cmidrule(lr){12-15}\cmidrule(lr){16-18}
$t$ & $q$ & $n$, $m$ & $|V|$ & $\alpha$ & time & cc & sc & $\alpha$ & time & cc & $\alpha$ & time & cc & sc & $\alpha$ & time & cc \\
\midrule
\endhead
\midrule
\multicolumn{18}{r}{Continued on next page} \\
\endfoot
\bottomrule
\endlastfoot
3 & 2 & $(1^{3})$, $(2,1,1)$ & 16 & \alphaopt{4} & 0 & 46 & 3 & \alphaopt{4} & 0 & 46 & \alphaopt{2} & 0 & 24 & 0 & \alphaopt{2} & 0 & 24 \\
2 & 2 & $(2,1)$, $(2,1)$ & 32 & \alphaopt{8} & 0 & 204 & 11 & \alphaopt{8} & 0 & 213 & \alphaopt{2} & 0 & 240 & 0 & \alphaopt{2} & 0 & 240 \\
4 & 2 & $(1^{4})$, $(2,1^{3})$ & 32 & \alphaopt{8} & 0 & 124 & 4 & \alphaopt{8} & 0 & 124 & \alphaopt{4} & 0 & 301 & 5 & \alphaopt{4} & 0 & 301 \\
5 & 2 & $(1^{5})$, $(2,1^{4})$ & 64 & \alphaopt{16} & 0 & 344 & 12 & \alphaopt{16} & 0 & 343 & \alphaopt{8} & 0 & 912 & 27 & \alphaopt{8} & 0 & 916 \\
2 & 2 & $(2,1)$, $(2,2)$ & 64 & \alphaopt{16} & 0 & 407 & 2 & \alphaopt{16} & 0 & 532 & \alphaopt{4} & 0 & 1461 & 9 & \alphaopt{4} & 0 & 1460 \\
6 & 2 & $(1^{6})$, $(2,1^{5})$ & 128 & \alphaopt{32} & 2 & 1008 & 12 & \alphaopt{32} & 2 & 1008 & \alphaopt{9} & 115 & 2736 & 3680 & \alphaopt{9} & 13 & 2705 \\
7 & 2 & $(1^{7})$, $(2,1^{6})$ & 256 & \alphaopt{64} & 20 & 2220 & 21 & \alphaopt{64} & 22 & 2195 & \alphanon{18} & 10840 & 7407 & 94 & \alphanon{18} & 10837 & 7305 \\
4 & 2 & $(2,1^{3})$, $(2,2,1,1)$ & 256 & \alphaopt{64} & 52 & 3980 & 34 & \alphaopt{64} & 48 & 4131 & \alphaopt{16} & 71 & 9792 & 81 & \alphaopt{16} & 64 & 9832 \\
\end{longtable}
\clearpage
\addtocounter{table}{-1}
\begin{longtable}{@{}r r l r !{\,\vrule\,} r r r r !{\,\vrule\,} r r r !{\,\vrule\,} r r r r !{\,\vrule\,} r r r @{}}
\caption{Sum-rank-metric graph instances: both cut families (Both) versus combinatorial cuts only (Comb.). Panel (b): $k=3,4$. Continued.}\\
\toprule
\multicolumn{4}{c}{} & \multicolumn{7}{c}{$k=3$} & \multicolumn{7}{c}{$k=4$} \\
\cmidrule(lr){5-11}\cmidrule(lr){12-18}
\multicolumn{4}{c}{} & \multicolumn{4}{c}{Both} & \multicolumn{3}{c}{Comb.} & \multicolumn{4}{c}{Both} & \multicolumn{3}{c}{Comb.} \\
\cmidrule(lr){5-8}\cmidrule(lr){9-11}\cmidrule(lr){12-15}\cmidrule(lr){16-18}
$t$ & $q$ & $n$, $m$ & $|V|$ & $\alpha$ & time & cc & sc & $\alpha$ & time & cc & $\alpha$ & time & cc & sc & $\alpha$ & time & cc \\
\midrule
\endfirsthead
\caption[]{Sum-rank-metric graph instances: both cut families (Both) versus combinatorial cuts only (Comb.). Panel (b), continued.}\\
\toprule
\multicolumn{4}{c}{} & \multicolumn{7}{c}{$k=3$} & \multicolumn{7}{c}{$k=4$} \\
\cmidrule(lr){5-11}\cmidrule(lr){12-18}
\multicolumn{4}{c}{} & \multicolumn{4}{c}{Both} & \multicolumn{3}{c}{Comb.} & \multicolumn{4}{c}{Both} & \multicolumn{3}{c}{Comb.} \\
\cmidrule(lr){5-8}\cmidrule(lr){9-11}\cmidrule(lr){12-15}\cmidrule(lr){16-18}
$t$ & $q$ & $n$, $m$ & $|V|$ & $\alpha$ & time & cc & sc & $\alpha$ & time & cc & $\alpha$ & time & cc & sc & $\alpha$ & time & cc \\
\midrule
\endhead
\midrule
\multicolumn{18}{r}{Continued on next page} \\
\endfoot
\bottomrule
\endlastfoot
3 & 2 & $(1^{3})$, $(2,1,1)$ & 16 & \alphaopt{1} & 0 & 0 & 0 & \alphaopt{1} & 0 & 0 & \alphaopt{1} & 0 & 0 & 0 & \alphaopt{1} & 0 & 0 \\
2 & 2 & $(2,1)$, $(2,1)$ & 32 & \alphaopt{1} & 0 & 0 & 0 & \alphaopt{1} & 0 & 0 & \alphaopt{1} & 0 & 0 & 0 & \alphaopt{1} & 0 & 0 \\
4 & 2 & $(1^{4})$, $(2,1^{3})$ & 32 & \alphaopt{2} & 0 & 48 & 0 & \alphaopt{2} & 0 & 48 & \alphaopt{1} & 0 & 0 & 0 & \alphaopt{1} & 0 & 0 \\
5 & 2 & $(1^{5})$, $(2,1^{4})$ & 64 & \alphaopt{2} & 0 & 1376 & 0 & \alphaopt{2} & 0 & 1376 & \alphaopt{2} & 0 & 96 & 0 & \alphaopt{2} & 0 & 96 \\
2 & 2 & $(2,1)$, $(2,2)$ & 64 & \alphaopt{1} & 0 & 0 & 0 & \alphaopt{1} & 0 & 0 & \alphaopt{1} & 0 & 0 & 0 & \alphaopt{1} & 0 & 0 \\
6 & 2 & $(1^{6})$, $(2,1^{5})$ & 128 & \alphaopt{4} & 7 & 4744 & 34 & \alphaopt{4} & 7 & 4727 & \alphaopt{2} & 6 & 4032 & 0 & \alphaopt{2} & 6 & 4032 \\
7 & 2 & $(1^{7})$, $(2,1^{6})$ & 256 & \alphaopt{8} & 87 & 14114 & 71 & \alphaopt{8} & 86 & 14125 & \alphaopt{4} & 146 & 23321 & 14 & \alphaopt{4} & 140 & 23321 \\
4 & 2 & $(2,1^{3})$, $(2,2,1,1)$ & 256 & \alphaopt{4} & 135 & 21615 & 59 & \alphaopt{4} & 119 & 21573 & \alphaopt{2} & 54 & 8064 & 0 & \alphaopt{2} & 46 & 8064 \\
\end{longtable}
\clearpage
\begin{longtable}{@{}l !{\,\vrule\,} r r r r !{\,\vrule\,} r r r !{\,\vrule\,} r r r r !{\,\vrule\,} r r r @{}}
\caption{Named graph instances: both cut families (Both) versus spectral cuts only (Spec.). Panel (a): $k=1,2$.}\label{tab:named-both-vs-spec}\\
\toprule
\multicolumn{1}{c}{} & \multicolumn{7}{c}{$k=1$} & \multicolumn{7}{c}{$k=2$} \\
\cmidrule(lr){2-8}\cmidrule(lr){9-15}
\multicolumn{1}{c}{} & \multicolumn{4}{c}{Both} & \multicolumn{3}{c}{Spec.} & \multicolumn{4}{c}{Both} & \multicolumn{3}{c}{Spec.} \\
\cmidrule(lr){2-5}\cmidrule(lr){6-8}\cmidrule(lr){9-12}\cmidrule(lr){13-15}
Instance & $\alpha$ & time & cc & sc & $\alpha$ & time & sc & $\alpha$ & time & cc & sc & $\alpha$ & time & sc \\
\midrule
\endfirsthead
\caption[]{Named graph instances: both cut families (Both) versus spectral cuts only (Spec.). Panel (a), continued.}\\
\toprule
\multicolumn{1}{c}{} & \multicolumn{7}{c}{$k=1$} & \multicolumn{7}{c}{$k=2$} \\
\cmidrule(lr){2-8}\cmidrule(lr){9-15}
\multicolumn{1}{c}{} & \multicolumn{4}{c}{Both} & \multicolumn{3}{c}{Spec.} & \multicolumn{4}{c}{Both} & \multicolumn{3}{c}{Spec.} \\
\cmidrule(lr){2-5}\cmidrule(lr){6-8}\cmidrule(lr){9-12}\cmidrule(lr){13-15}
Instance & $\alpha$ & time & cc & sc & $\alpha$ & time & sc & $\alpha$ & time & cc & sc & $\alpha$ & time & sc \\
\midrule
\endhead
\midrule
\multicolumn{15}{r}{Continued on next page} \\
\endfoot
\bottomrule
\endlastfoot
Tetrahedral & \alphaopt{1} & 0 & 0 & 0 & \alphaopt{1} & 0 & 0 & \alphaopt{1} & 0 & 0 & 0 & \alphaopt{1} & 0 & 0 \\
Thomsen & \alphaopt{3} & 0 & 0 & 0 & \alphaopt{3} & 0 & 0 & \alphaopt{1} & 0 & 0 & 0 & \alphaopt{1} & 0 & 0 \\
Octahedral & \alphaopt{2} & 0 & 0 & 0 & \alphaopt{2} & 0 & 0 & \alphaopt{1} & 0 & 0 & 0 & \alphaopt{1} & 0 & 0 \\
MoserSpindle & \alphaopt{2} & 0 & 11 & 0 & \alphaopt{2} & 0 & 0 & \alphaopt{1} & 0 & 0 & 0 & \alphaopt{1} & 0 & 0 \\
Hexahedral & \alphaopt{4} & 0 & 12 & 0 & \alphaopt{4} & 0 & 0 & \alphaopt{2} & 0 & 0 & 0 & \alphaopt{2} & 0 & 0 \\
Wagner & \alphaopt{3} & 0 & 12 & 0 & \alphaopt{3} & 0 & 0 & \alphaopt{1} & 0 & 0 & 0 & \alphaopt{1} & 0 & 0 \\
Petersen & \alphaopt{4} & 0 & 15 & 0 & \alphaopt{4} & 0 & 0 & \alphaopt{1} & 0 & 0 & 0 & \alphaopt{1} & 0 & 0 \\
KrackhardtKite & \alphaopt{4} & 0 & 20 & 1 & \alphaopt{4} & 0 & 4 & \alphaopt{2} & 0 & 22 & 0 & \alphaopt{2} & 0 & 0 \\
\end{longtable}
\clearpage
\addtocounter{table}{-1}
\begin{longtable}{@{}l !{\,\vrule\,} r r r r !{\,\vrule\,} r r r !{\,\vrule\,} r r r r !{\,\vrule\,} r r r @{}}
\caption{Named graph instances: both cut families (Both) versus spectral cuts only (Spec.). Panel (b): $k=3,4$. Continued.}\\
\toprule
\multicolumn{1}{c}{} & \multicolumn{7}{c}{$k=3$} & \multicolumn{7}{c}{$k=4$} \\
\cmidrule(lr){2-8}\cmidrule(lr){9-15}
\multicolumn{1}{c}{} & \multicolumn{4}{c}{Both} & \multicolumn{3}{c}{Spec.} & \multicolumn{4}{c}{Both} & \multicolumn{3}{c}{Spec.} \\
\cmidrule(lr){2-5}\cmidrule(lr){6-8}\cmidrule(lr){9-12}\cmidrule(lr){13-15}
Instance & $\alpha$ & time & cc & sc & $\alpha$ & time & sc & $\alpha$ & time & cc & sc & $\alpha$ & time & sc \\
\midrule
\endfirsthead
\caption[]{Named graph instances: both cut families (Both) versus spectral cuts only (Spec.). Panel (b), continued.}\\
\toprule
\multicolumn{1}{c}{} & \multicolumn{7}{c}{$k=3$} & \multicolumn{7}{c}{$k=4$} \\
\cmidrule(lr){2-8}\cmidrule(lr){9-15}
\multicolumn{1}{c}{} & \multicolumn{4}{c}{Both} & \multicolumn{3}{c}{Spec.} & \multicolumn{4}{c}{Both} & \multicolumn{3}{c}{Spec.} \\
\cmidrule(lr){2-5}\cmidrule(lr){6-8}\cmidrule(lr){9-12}\cmidrule(lr){13-15}
Instance & $\alpha$ & time & cc & sc & $\alpha$ & time & sc & $\alpha$ & time & cc & sc & $\alpha$ & time & sc \\
\midrule
\endhead
\midrule
\multicolumn{15}{r}{Continued on next page} \\
\endfoot
\bottomrule
\endlastfoot
Tetrahedral & \alphaopt{1} & 0 & 0 & 0 & \alphaopt{1} & 0 & 0 & \alphaopt{1} & 0 & 0 & 0 & \alphaopt{1} & 0 & 0 \\
Thomsen & \alphaopt{1} & 0 & 0 & 0 & \alphaopt{1} & 0 & 0 & \alphaopt{1} & 0 & 0 & 0 & \alphaopt{1} & 0 & 0 \\
Octahedral & \alphaopt{1} & 0 & 0 & 0 & \alphaopt{1} & 0 & 0 & \alphaopt{1} & 0 & 0 & 0 & \alphaopt{1} & 0 & 0 \\
MoserSpindle & \alphaopt{1} & 0 & 0 & 0 & \alphaopt{1} & 0 & 0 & \alphaopt{1} & 0 & 0 & 0 & \alphaopt{1} & 0 & 0 \\
Hexahedral & \alphaopt{1} & 0 & 0 & 0 & \alphaopt{1} & 0 & 0 & \alphaopt{1} & 0 & 0 & 0 & \alphaopt{1} & 0 & 0 \\
Wagner & \alphaopt{1} & 0 & 0 & 0 & \alphaopt{1} & 0 & 0 & \alphaopt{1} & 0 & 0 & 0 & \alphaopt{1} & 0 & 0 \\
Petersen & \alphaopt{1} & 0 & 0 & 0 & \alphaopt{1} & 0 & 0 & \alphaopt{1} & 0 & 0 & 0 & \alphaopt{1} & 0 & 0 \\
KrackhardtKite & \alphaopt{2} & 0 & 10 & 0 & \alphaopt{2} & 0 & 0 & \alphaopt{1} & 0 & 0 & 0 & \alphaopt{1} & 0 & 0 \\
\end{longtable}
\clearpage
\begin{longtable}{@{}l !{\,\vrule\,} r r r r !{\,\vrule\,} r r r !{\,\vrule\,} r r r r !{\,\vrule\,} r r r @{}}
\caption{Coding-theoretic instances: both cut families (Both) versus spectral cuts only (Spec.). Panel (a): $k=1,2$.}\label{tab:numeric-both-vs-spec}\\
\toprule
\multicolumn{1}{c}{} & \multicolumn{7}{c}{$k=1$} & \multicolumn{7}{c}{$k=2$} \\
\cmidrule(lr){2-8}\cmidrule(lr){9-15}
\multicolumn{1}{c}{} & \multicolumn{4}{c}{Both} & \multicolumn{3}{c}{Spec.} & \multicolumn{4}{c}{Both} & \multicolumn{3}{c}{Spec.} \\
\cmidrule(lr){2-5}\cmidrule(lr){6-8}\cmidrule(lr){9-12}\cmidrule(lr){13-15}
Instance & $\alpha$ & time & cc & sc & $\alpha$ & time & sc & $\alpha$ & time & cc & sc & $\alpha$ & time & sc \\
\midrule
\endfirsthead
\caption[]{Coding-theoretic instances: both cut families (Both) versus spectral cuts only (Spec.). Panel (a), continued.}\\
\toprule
\multicolumn{1}{c}{} & \multicolumn{7}{c}{$k=1$} & \multicolumn{7}{c}{$k=2$} \\
\cmidrule(lr){2-8}\cmidrule(lr){9-15}
\multicolumn{1}{c}{} & \multicolumn{4}{c}{Both} & \multicolumn{3}{c}{Spec.} & \multicolumn{4}{c}{Both} & \multicolumn{3}{c}{Spec.} \\
\cmidrule(lr){2-5}\cmidrule(lr){6-8}\cmidrule(lr){9-12}\cmidrule(lr){13-15}
Instance & $\alpha$ & time & cc & sc & $\alpha$ & time & sc & $\alpha$ & time & cc & sc & $\alpha$ & time & sc \\
\midrule
\endhead
\midrule
\multicolumn{15}{r}{Continued on next page} \\
\endfoot
\bottomrule
\endlastfoot
1tc.8 & \alphaopt{4} & 0 & 6 & 0 & \alphaopt{4} & 0 & 0 & \alphaopt{4} & 0 & 6 & 0 & \alphaopt{4} & 0 & 0 \\
1tc.16 & \alphaopt{8} & 0 & 35 & 5 & \alphaopt{8} & 0 & 3 & \alphaopt{5} & 0 & 27 & 0 & \alphaopt{5} & 0 & 0 \\
1tc.32 & \alphaopt{12} & 0 & 98 & 3 & \alphaopt{12} & 0 & 197 & \alphaopt{8} & 0 & 125 & 4 & \alphaopt{8} & 0 & 44 \\
1dc.64 & \alphaopt{10} & 0 & 659 & 30 & \alphaopt{10} & 235 & 24651 & \alphaopt{4} & 0 & 1398 & 2 & \alphaopt{4} & 0 & 72 \\
1et.64 & \alphaopt{18} & 0 & 635 & 30 & \alphaopt{18} & 2 & 827 & \alphaopt{8} & 0 & 452 & 6 & \alphaopt{8} & 12 & 5785 \\
1tc.64 & \alphaopt{20} & 0 & 478 & 13 & \alphaopt{20} & 180 & 16292 & \alphaopt{12} & 0 & 467 & 18 & \alphaopt{12} & 5 & 3217 \\
1dc.128 & \alphaopt{16} & 4 & 2017 & 50 & \alphanon{16} & 14582 & 166669 & \alphaopt{5} & 7 & 5237 & 24 & \alphaopt{5} & 145 & 17495 \\
1et.128 & \alphaopt{28} & 3 & 1639 & 26 & \alphanon{28} & 10806 & 113859 & \alphaopt{14} & 3 & 1717 & 20 & \alphanon{14} & 10801 & 120016 \\
\end{longtable}
\clearpage
\addtocounter{table}{-1}
\begin{longtable}{@{}l !{\,\vrule\,} r r r r !{\,\vrule\,} r r r !{\,\vrule\,} r r r r !{\,\vrule\,} r r r @{}}
\caption{Coding-theoretic instances: both cut families (Both) versus spectral cuts only (Spec.). Panel (b): $k=3,4$. Continued.}\\
\toprule
\multicolumn{1}{c}{} & \multicolumn{7}{c}{$k=3$} & \multicolumn{7}{c}{$k=4$} \\
\cmidrule(lr){2-8}\cmidrule(lr){9-15}
\multicolumn{1}{c}{} & \multicolumn{4}{c}{Both} & \multicolumn{3}{c}{Spec.} & \multicolumn{4}{c}{Both} & \multicolumn{3}{c}{Spec.} \\
\cmidrule(lr){2-5}\cmidrule(lr){6-8}\cmidrule(lr){9-12}\cmidrule(lr){13-15}
Instance & $\alpha$ & time & cc & sc & $\alpha$ & time & sc & $\alpha$ & time & cc & sc & $\alpha$ & time & sc \\
\midrule
\endfirsthead
\caption[]{Coding-theoretic instances: both cut families (Both) versus spectral cuts only (Spec.). Panel (b), continued.}\\
\toprule
\multicolumn{1}{c}{} & \multicolumn{7}{c}{$k=3$} & \multicolumn{7}{c}{$k=4$} \\
\cmidrule(lr){2-8}\cmidrule(lr){9-15}
\multicolumn{1}{c}{} & \multicolumn{4}{c}{Both} & \multicolumn{3}{c}{Spec.} & \multicolumn{4}{c}{Both} & \multicolumn{3}{c}{Spec.} \\
\cmidrule(lr){2-5}\cmidrule(lr){6-8}\cmidrule(lr){9-12}\cmidrule(lr){13-15}
Instance & $\alpha$ & time & cc & sc & $\alpha$ & time & sc & $\alpha$ & time & cc & sc & $\alpha$ & time & sc \\
\midrule
\endhead
\midrule
\multicolumn{15}{r}{Continued on next page} \\
\endfoot
\bottomrule
\endlastfoot
1tc.8 & \alphaopt{4} & 0 & 6 & 0 & \alphaopt{4} & 0 & 0 & \alphaopt{4} & 0 & 6 & 0 & \alphaopt{4} & 0 & 0 \\
1tc.16 & \alphaopt{5} & 0 & 27 & 0 & \alphaopt{5} & 0 & 0 & \alphaopt{5} & 0 & 27 & 0 & \alphaopt{5} & 0 & 0 \\
1tc.32 & \alphaopt{6} & 0 & 110 & 0 & \alphaopt{6} & 0 & 0 & \alphaopt{6} & 0 & 110 & 0 & \alphaopt{6} & 0 & 0 \\
1dc.64 & \alphaopt{2} & 0 & 462 & 0 & \alphaopt{2} & 0 & 0 & \alphaopt{2} & 0 & 42 & 0 & \alphaopt{2} & 0 & 0 \\
1et.64 & \alphaopt{7} & 0 & 430 & 0 & \alphaopt{7} & 0 & 0 & \alphaopt{7} & 0 & 430 & 0 & \alphaopt{7} & 0 & 0 \\
1tc.64 & \alphaopt{10} & 0 & 460 & 6 & \alphaopt{10} & 4 & 3109 & \alphaopt{8} & 0 & 435 & 1 & \alphaopt{8} & 0 & 21 \\
1dc.128 & \alphaopt{2} & 6 & 4032 & 0 & \alphaopt{2} & 0 & 0 & \alphaopt{2} & 1 & 812 & 0 & \alphaopt{2} & 0 & 0 \\
1et.128 & \alphaopt{8} & 3 & 1652 & 0 & \alphaopt{8} & 0 & 0 & \alphaopt{8} & 3 & 1652 & 0 & \alphaopt{8} & 0 & 0 \\
\end{longtable}
\clearpage
\begin{longtable}{@{}r r l r !{\,\vrule\,} r r r r !{\,\vrule\,} r r r !{\,\vrule\,} r r r r !{\,\vrule\,} r r r @{}}
\caption{Sum-rank-metric graph instances: both cut families (Both) versus spectral cuts only (Spec.). Panel (a): $k=1,2$.}\label{tab:graph-both-vs-spec}\\
\toprule
\multicolumn{4}{c}{} & \multicolumn{7}{c}{$k=1$} & \multicolumn{7}{c}{$k=2$} \\
\cmidrule(lr){5-11}\cmidrule(lr){12-18}
\multicolumn{4}{c}{} & \multicolumn{4}{c}{Both} & \multicolumn{3}{c}{Spec.} & \multicolumn{4}{c}{Both} & \multicolumn{3}{c}{Spec.} \\
\cmidrule(lr){5-8}\cmidrule(lr){9-11}\cmidrule(lr){12-15}\cmidrule(lr){16-18}
$t$ & $q$ & $n$, $m$ & $|V|$ & $\alpha$ & time & cc & sc & $\alpha$ & time & sc & $\alpha$ & time & cc & sc & $\alpha$ & time & sc \\
\midrule
\endfirsthead
\caption[]{Sum-rank-metric graph instances: both cut families (Both) versus spectral cuts only (Spec.). Panel (a), continued.}\\
\toprule
\multicolumn{4}{c}{} & \multicolumn{7}{c}{$k=1$} & \multicolumn{7}{c}{$k=2$} \\
\cmidrule(lr){5-11}\cmidrule(lr){12-18}
\multicolumn{4}{c}{} & \multicolumn{4}{c}{Both} & \multicolumn{3}{c}{Spec.} & \multicolumn{4}{c}{Both} & \multicolumn{3}{c}{Spec.} \\
\cmidrule(lr){5-8}\cmidrule(lr){9-11}\cmidrule(lr){12-15}\cmidrule(lr){16-18}
$t$ & $q$ & $n$, $m$ & $|V|$ & $\alpha$ & time & cc & sc & $\alpha$ & time & sc & $\alpha$ & time & cc & sc & $\alpha$ & time & sc \\
\midrule
\endhead
\midrule
\multicolumn{18}{r}{Continued on next page} \\
\endfoot
\bottomrule
\endlastfoot
3 & 2 & $(1^{3})$, $(2,1,1)$ & 16 & \alphaopt{4} & 0 & 46 & 3 & \alphaopt{4} & 0 & 11 & \alphaopt{2} & 0 & 24 & 0 & \alphaopt{2} & 0 & 0 \\
2 & 2 & $(2,1)$, $(2,1)$ & 32 & \alphaopt{8} & 0 & 204 & 11 & \alphaopt{8} & 0 & 217 & \alphaopt{2} & 0 & 240 & 0 & \alphaopt{2} & 0 & 0 \\
4 & 2 & $(1^{4})$, $(2,1^{3})$ & 32 & \alphaopt{8} & 0 & 124 & 4 & \alphaopt{8} & 0 & 323 & \alphaopt{4} & 0 & 301 & 5 & \alphaopt{4} & 0 & 11 \\
5 & 2 & $(1^{5})$, $(2,1^{4})$ & 64 & \alphaopt{16} & 0 & 344 & 12 & \alphaopt{16} & 289 & 25548 & \alphaopt{8} & 0 & 912 & 27 & \alphaopt{8} & 235 & 11122 \\
2 & 2 & $(2,1)$, $(2,2)$ & 64 & \alphaopt{16} & 0 & 407 & 2 & \alphaopt{16} & 0 & 11 & \alphaopt{4} & 0 & 1461 & 9 & \alphaopt{4} & 0 & 9 \\
6 & 2 & $(1^{6})$, $(2,1^{5})$ & 128 & \alphaopt{32} & 2 & 1008 & 12 & \alphanon{32} & 10801 & 105624 & \alphaopt{9} & 115 & 2736 & 3680 & \alphanon{9} & 10801 & 140090 \\
7 & 2 & $(1^{7})$, $(2,1^{6})$ & 256 & \alphaopt{64} & 20 & 2220 & 21 & \alphanon{64} & 10825 & 1035 & \alphanon{18} & 10840 & 7407 & 94 & \alphanon{18} & 13657 & 157374 \\
4 & 2 & $(2,1^{3})$, $(2,2,1,1)$ & 256 & \alphaopt{64} & 52 & 3980 & 34 & \alphanon{56} & 11538 & 45593 & \alphaopt{16} & 71 & 9792 & 81 & \alphanon{16} & 10803 & 63546 \\
\end{longtable}
\clearpage
\addtocounter{table}{-1}
\begin{longtable}{@{}r r l r !{\,\vrule\,} r r r r !{\,\vrule\,} r r r !{\,\vrule\,} r r r r !{\,\vrule\,} r r r @{}}
\caption{Sum-rank-metric graph instances: both cut families (Both) versus spectral cuts only (Spec.). Panel (b): $k=3,4$. Continued.}\\
\toprule
\multicolumn{4}{c}{} & \multicolumn{7}{c}{$k=3$} & \multicolumn{7}{c}{$k=4$} \\
\cmidrule(lr){5-11}\cmidrule(lr){12-18}
\multicolumn{4}{c}{} & \multicolumn{4}{c}{Both} & \multicolumn{3}{c}{Spec.} & \multicolumn{4}{c}{Both} & \multicolumn{3}{c}{Spec.} \\
\cmidrule(lr){5-8}\cmidrule(lr){9-11}\cmidrule(lr){12-15}\cmidrule(lr){16-18}
$t$ & $q$ & $n$, $m$ & $|V|$ & $\alpha$ & time & cc & sc & $\alpha$ & time & sc & $\alpha$ & time & cc & sc & $\alpha$ & time & sc \\
\midrule
\endfirsthead
\caption[]{Sum-rank-metric graph instances: both cut families (Both) versus spectral cuts only (Spec.). Panel (b), continued.}\\
\toprule
\multicolumn{4}{c}{} & \multicolumn{7}{c}{$k=3$} & \multicolumn{7}{c}{$k=4$} \\
\cmidrule(lr){5-11}\cmidrule(lr){12-18}
\multicolumn{4}{c}{} & \multicolumn{4}{c}{Both} & \multicolumn{3}{c}{Spec.} & \multicolumn{4}{c}{Both} & \multicolumn{3}{c}{Spec.} \\
\cmidrule(lr){5-8}\cmidrule(lr){9-11}\cmidrule(lr){12-15}\cmidrule(lr){16-18}
$t$ & $q$ & $n$, $m$ & $|V|$ & $\alpha$ & time & cc & sc & $\alpha$ & time & sc & $\alpha$ & time & cc & sc & $\alpha$ & time & sc \\
\midrule
\endhead
\midrule
\multicolumn{18}{r}{Continued on next page} \\
\endfoot
\bottomrule
\endlastfoot
3 & 2 & $(1^{3})$, $(2,1,1)$ & 16 & \alphaopt{1} & 0 & 0 & 0 & \alphaopt{1} & 0 & 0 & \alphaopt{1} & 0 & 0 & 0 & \alphaopt{1} & 0 & 0 \\
2 & 2 & $(2,1)$, $(2,1)$ & 32 & \alphaopt{1} & 0 & 0 & 0 & \alphaopt{1} & 0 & 0 & \alphaopt{1} & 0 & 0 & 0 & \alphaopt{1} & 0 & 0 \\
4 & 2 & $(1^{4})$, $(2,1^{3})$ & 32 & \alphaopt{2} & 0 & 48 & 0 & \alphaopt{2} & 0 & 0 & \alphaopt{1} & 0 & 0 & 0 & \alphaopt{1} & 0 & 0 \\
5 & 2 & $(1^{5})$, $(2,1^{4})$ & 64 & \alphaopt{2} & 0 & 1376 & 0 & \alphaopt{2} & 0 & 0 & \alphaopt{2} & 0 & 96 & 0 & \alphaopt{2} & 0 & 0 \\
2 & 2 & $(2,1)$, $(2,2)$ & 64 & \alphaopt{1} & 0 & 0 & 0 & \alphaopt{1} & 0 & 0 & \alphaopt{1} & 0 & 0 & 0 & \alphaopt{1} & 0 & 0 \\
6 & 2 & $(1^{6})$, $(2,1^{5})$ & 128 & \alphaopt{4} & 7 & 4744 & 34 & \alphaopt{4} & 2 & 1866 & \alphaopt{2} & 6 & 4032 & 0 & \alphaopt{2} & 0 & 0 \\
7 & 2 & $(1^{7})$, $(2,1^{6})$ & 256 & \alphaopt{8} & 87 & 14114 & 71 & \alphaopt{8} & 415 & 17026 & \alphaopt{4} & 146 & 23321 & 14 & \alphaopt{4} & 0 & 29 \\
4 & 2 & $(2,1^{3})$, $(2,2,1,1)$ & 256 & \alphaopt{4} & 135 & 21615 & 59 & \alphaopt{4} & 37 & 5408 & \alphaopt{2} & 54 & 8064 & 0 & \alphaopt{2} & 0 & 0 \\
\end{longtable}

\end{landscape}
\endgroup

\end{document}